\documentclass[12pt,a4paper]{amsart}
\usepackage[T1]{fontenc}
\usepackage[utf8]{inputenc}
\usepackage{lmodern}
\usepackage[english]{babel}
\usepackage[margin=1.2in]{geometry}
\usepackage{amsmath,amssymb,amsthm,mathtools,stmaryrd,mathrsfs}
\usepackage{microtype}
\usepackage{indentfirst}
\usepackage{needspace}
\usepackage{enumitem}
\usepackage[numbers,sort&compress]{natbib}
\usepackage[colorlinks=true,linkcolor=black,citecolor=black,urlcolor=black,
  bookmarksnumbered=true,hypertexnames=false]{hyperref}
\usepackage{bookmark}
\allowdisplaybreaks[2]
\numberwithin{equation}{section}
\newtheorem{theorem}{Theorem}[section]
\newtheorem{proposition}[theorem]{Proposition}
\newtheorem{lemma}[theorem]{Lemma}
\newtheorem{corollary}[theorem]{Corollary}
\theoremstyle{definition}
\newtheorem{definition}[theorem]{Definition}

\theoremstyle{remark}
\newtheorem{remark}[theorem]{Remark}

\newcommand{\Zt}{\mathbb Z_2}
\newcommand{\R}{\mathbb R}

\newcommand{\Mass}{\mathbf M}
\newcommand{\Flat}{\mathcal F_r}
\newcommand{\dr}{\partial_r}
\newcommand{\II}{\mathbf I}
\newcommand{\IC}{\mathcal I}
\newcommand{\ZC}{\mathcal Z}
\newcommand{\VC}{\mathcal V}
\newcommand{\FH}{\mathbf F}
\newcommand{\AMset}{\mathscr A}
\newcommand{\Alm}{\mathcal A}
\newcommand{\MC}{\llbracket M\rrbracket}
\newcommand{\res}{\mathbin{\llcorner}}
\newcommand{\HH}{\mathcal H}
\newcommand{\spt}{\operatorname{spt}}
\newcommand{\Int}{\operatorname{Int}}
\newcommand{\Vol}{\operatorname{Vol}}
\newcommand{\diam}{\operatorname{diam}}
\newcommand{\dist}{\operatorname{dist}}
\newcommand{\Lip}{\operatorname{Lip}}
\newcommand{\diver}{\operatorname{div}}
\newcommand{\SFF}{\mathrm{II}}
\newcommand{\dV}{\mathbf d_V}
\newcommand{\eps}{\varepsilon}
\newcommand{\An}{\operatorname{An}}

\hypersetup{
  pdftitle={Free-Boundary Minimal Varifolds II: Relative Filling Bounds},
  pdfauthor={Tengyao Zhang and Zhifei Zhu},
  pdfsubject={Relative filling functions, free-boundary min--max, mass bounds}
}
\title{Free-Boundary Minimal Varifolds II:\\ Relative Filling Bounds}
\author{Tengyao Zhang}
\address{Yau Mathematical Sciences Center, Tsinghua University, Beijing, 100084, P.R. China}
\email{zhangtengyao25@mails.tsinghua.edu.cn}

\author{Zhifei Zhu}
\address{Yau Mathematical Sciences Center, Tsinghua University, Beijing, 100084, P.R. China}
\email{zhifeizhu@tsinghua.edu.cn}
\date{}
\begin{document}
\begin{abstract}
Let $M^n$ be a compact connected Riemannian manifold with nonempty smooth boundary, and let $2\le k<n$. Assuming $H_j(M,\partial M;\Zt)=0$ for $1\le j<k$, we construct a nonzero free-boundary stationary integral $k$-varifold whose mass is bounded in terms of the relative filling functions of $M$ and the maximal distance to the boundary. The varifold is almost minimizing modulo two in small relative annuli. By the first paper of this series, strict trace conditions on the second fundamental form of the boundary therefore yield zero boundary mass; for hypersurfaces with $3\le n\le7$, strict mean convexity gives smooth, properly embedded free-boundary minimal hypersurfaces. In dimension four we bound the mass of a stationary two-varifold in terms of Ricci curvature, volume and diameter bounds for a smooth metric on the double of $M$ and its bilipschitz comparison constant with the given metric on each half.
\end{abstract}
\maketitle

\setcounter{tocdepth}{2}
\begingroup
\small
\tableofcontents
\endgroup
\clearpage

\section{Introduction}
\label{sec:introduction}

Let $(M^n,g)$ be a smooth compact connected Riemannian manifold of
dimension $n\ge3$ with nonempty smooth boundary, and let
$2\le k\le n-1$. A $k$-varifold in $M$ is said to be \emph{stationary
with free boundary} if its first variation vanishes on every smooth
vector field tangent to $\partial M$. For a smooth properly embedded
submanifold, this condition is equivalent to minimality together with
orthogonal contact along the boundary. This paper is concerned with the
existence of nonzero $k$-varifolds of this kind and with upper bounds
for their mass in terms of the geometry of $M$.

On a closed manifold, the min--max theory of Almgren and Pitts produces
a nonzero stationary integral $k$-varifold for every $k$
\cite{Almgren1965,Pitts1981}. When the homology of the manifold
vanishes in degrees one through $k-1$, Nabutovsky and Rotman bounded
the least mass of such a varifold in terms of the diameter and the
homological filling functions \cite{NR2006}. The value of a filling
function at $L$ is the sharp upper bound on the volume needed to fill
a cycle of mass at most $L$. The study of such bounds goes back to
Gromov \cite{Gromov1983}.
Nabutovsky, Rotman, and Sabourau developed related sweepout and
homological-waist constructions \cite{NRS2020}.

We extend the estimate of Nabutovsky and Rotman to manifolds with
boundary. The filling recursion and the fundamental-class obstruction
are carried out in the relative cycle spaces of $(M,\partial M)$, and
we establish the compactness needed for the limiting varifold. Joining
vertices to nearest boundary points replaces the diameter in the first
step of the recursion by the largest distance to the boundary,
\begin{equation}\label{eq:boundarydistance}
 D=\sup_{x\in M}\dist_g(x,\partial M)>0,
\end{equation}
where distance is measured along curves in $M$; the factorial
coefficient remains that of \cite{NR2006}. If
$H_j(M,\partial M;\Zt)=0$ for $1\le j<k$, we obtain a nonzero integral
$k$-varifold, stationary with free boundary and almost minimizing
modulo two in small relative annuli, whose mass is bounded in terms of
$D$ and the relative filling functions of $M$
(Theorem~\ref{thm:mass}). In dimension four, a two-sided Ricci bound, a
lower volume bound and an upper diameter bound on the smooth double,
together with a uniform bilipschitz comparison on each half, give a
uniform bound for the first relative filling function at each length
and hence for the mass of a stationary two-varifold
(Theorem~\ref{thm:relative-wz}).

In dimension three, suppose $H_1(M,\partial M;\Zt)=0$ and every
relative one-cycle of length at most $2D$ has a relative filling of
area at most $A$. Theorem~\ref{thm:mass} then gives a two-varifold
of mass at most $12A$. If the
boundary is strictly mean convex, this gives a smooth, compact,
connected, properly embedded free-boundary minimal surface of area at
most $12A$, possibly with empty boundary, by Part~I
\cite[Corollary~1.6]{PartI}.

Throughout the paper, $(M^n,g)$ is as above. Chains and homology groups
have coefficients in $\Zt$ unless otherwise stated. Stationarity with
free boundary and almost minimization in small relative annuli are
defined in Section~\ref{sec:notation}.

\subsection{Main results}

The mass bound is stated in terms of the following functions.

\begin{definition}
\label{def:filling}
For $1\le j\le n-1$ and $v\ge0$, define
\begin{equation}\label{eq:filling}
 \FH_j(v)=
 \sup_{\substack{z\in\ZC_j(M,\partial M;\Zt)\\\Mass(z)\le v}}
 \inf_{\substack{Q\in\IC_{j+1}(M,\partial M;\Zt)\\\dr Q=z}}
 \Mass(Q).
\end{equation}
Here $\ZC_j(M,\partial M;\Zt)$ is the space of relative $j$-cycles, $\IC_{j+1}(M,\partial M;\Zt)$ is the space of relative $(j+1)$-chains, $\dr$ is the relative boundary operator, and $\Mass$ is relative mass; these objects are defined in Section~\ref{sec:notation}. The infimum over the empty set is $+\infty$.
\end{definition}

We restrict our main results to $2\le k\le n-1$; the one-dimensional
case is discussed in Subsection~\ref{sec:reference-comparison}. The
first main result concerns mass alone and imposes no condition on the
second fundamental form of the boundary.

\begin{theorem}
\label{thm:mass}
Let $2\le k\le n-1$, and assume
\begin{equation}\label{eq:homology}
 H_j(M,\partial M;\Zt)=0
 \quad(1\le j\le k-1).
\end{equation}
Define
\begin{equation}\label{eq:budgets}
 b_1=\FH_1(2D),\quad
 b_j=\FH_j\bigl((j+1)b_{j-1}\bigr)
 \quad(2\le j\le k-1).
\end{equation}
Every filling function occurring here is finite and right-continuous, and all the numbers $b_j$ are positive. There is a nonzero integral $k$-varifold $V$ in $M$, stationary with free boundary and almost minimizing modulo two in small relative annuli, such that
\begin{equation}\label{eq:mainmass}
 0<\|V\|(M)\le\frac{(n+1)!}{k!}\,b_{k-1}.
\end{equation}
\end{theorem}

The argument $2D$ comes from joining the two endpoints of an edge to
$\partial M$, and the later factors $j+1$ count the faces of a simplex.
The relative fundamental class obstructs the completion of all these
fillings. The construction uses at most $n-k$ parameters and requires
no vanishing of $H_k(M,\partial M;\Zt)$ or $H_{k+1}(M,\partial M;\Zt)$;
see Section~\ref{sec:mass}. It follows Almgren
\cite{Almgren1962,Almgren1965}, Pitts \cite{Pitts1981}, and Li--Zhou
\cite{LZ2021}, with the relative interpolation of
Liokumovich--Marques--Neves \cite{LMN2018} for cycles of general
dimension. In particular, Theorem~\ref{thm:mass} does not need the
restriction on the ambient dimension used for smooth hypersurface
regularity.

Let $V$ be the varifold of Theorem~\ref{thm:mass} and write
$\SFF_{\partial M}(v,w)=g(\nabla_v\nu,w)$ for the second fundamental
form with respect to the outward unit normal $\nu$. If at each $x\in\partial M$
the traces of $\SFF_{\partial M}$ on all $k$-planes in $T_x\partial M$
have the same strict sign, then $\|V\|(\partial M)=0$ by
\cite[Corollary~1.5]{PartI}. The sign may vary between boundary
components.
For a smooth bounded Euclidean domain $\Omega$, even without a sign
condition, $\|V\|(\Omega)\ge c_{\Omega,k}\|V\|(\overline\Omega)$ for
some $c_{\Omega,k}>0$ by \cite[Corollary~1.7]{PartI}.

For $k=n-1$, set $H_{\partial M}=\operatorname{tr}_{T\partial M}
\SFF_{\partial M}$. If $3\le n\le7$ and $H_{\partial M}>0$, then
$V$ is represented by finitely many smooth, compact, properly embedded free-boundary
minimal hypersurfaces with positive integer multiplicities
\cite[Corollary~1.6]{PartI}. Their boundaries may be empty, and their
total area counted with multiplicity is at most $n(n+1)b_{n-2}$ by
\eqref{eq:mainmass}.

In dimension four, geometric bounds on the double yield estimates for
the relative filling functions. For closed four-manifolds, Wu and Zhu \cite{WZ2020}
bound the integer filling function $HF_1$ under two-sided Ricci bounds,
noncollapse, bounded diameter, and vanishing integer $H_1$. We prove a
mod-two version of their estimate for bounding cycles, without any
assumption on absolute homology (Proposition~\ref{wz:prop:closed}), and
we transfer it to manifolds with boundary through the double
(Proposition~\ref{wz:prop:transfer}). Combining the two estimates
with Theorem~\ref{thm:mass} gives the following result.

\Needspace{5\baselineskip}
\begin{theorem}\label{thm:relative-wz}
For every $v_0,d_0>0$ there are constants $a=a(v_0,d_0)>0$ and
$b=b(v_0,d_0)>0$ with the following property.
Let $(M^4,g)$ be smooth, compact and connected, with nonempty smooth
boundary and $H_1(M,\partial M;\Zt)=0$. Suppose its smooth double
$X=DM$ carries a smooth metric $h$ such that
\begin{equation}\label{wz:eq:double-hyp}
 |\operatorname{Ric}_h|\le3,\quad
 \Vol_h(X)\ge v_0,\quad \diam_h(X)\le d_0,
 \quad
 \Lambda^{-2}g\le j_\pm^*h\le\Lambda^2g
\end{equation}
for some $\Lambda\ge1$, where $j_\pm$ identify $M$ with the two
halves. Then
\begin{equation}\label{wz:eq:relative-main}
 \FH_1(L)\le \Lambda^3aL+\frac{\Lambda^2}{2}b
 \quad(L\ge0).
\end{equation}
Consequently there is a nonzero integral two-varifold $V$, stationary
with free boundary and almost minimizing modulo two in small relative
annuli, satisfying
\begin{equation}\label{wz:eq:mass-main}
 0<\|V\|(M)\le120\Lambda^3aD+30\Lambda^2b,
\end{equation}
where $D$ is defined in \eqref{eq:boundarydistance}. If at each
boundary point the two-plane traces of $\SFF_{\partial M}$ have the
same strict sign, then $\|V\|(\partial M)=0$.
\end{theorem}

The metric on the double need not be invariant under the interchange
of its halves, and no assumption is imposed on the absolute first
homology of the double. Proposition~\ref{wz:prop:obstruction} shows
that two-sided Ricci bounds, a lower volume bound, and an upper
diameter bound on $M$ do not give a uniform upper bound for $\FH_1(2D)$.

The same method gives more than is needed for Theorem~\ref{thm:relative-wz}.
The closed body--neck argument also gives an affine mod-two $HF_2$
estimate for bounding cycles and, on oriented closed four-manifolds, a
linear integer estimate (Proposition~\ref{hf:prop:closed2});
top-dimensional filling gives the $HF_3$ bound \eqref{hf:eq:closed3}.
These closed estimates are proved in Section~\ref{sec:closed-fill}.
Section~\ref{sec:relative-wz} proves a transfer principle in every
degree and applies it to the three closed estimates; a bound for all
cycles of a given degree requires vanishing of the relative homology
group in that degree. In particular, three-dimensional hypersurfaces in
a four-manifold are controlled by $\FH_1$ and $\FH_2$:
Corollary~\ref{hf:cor:hypersurface} gives the mass estimate and, under
strict mean convexity, smooth proper embeddedness. Product-collar and
puncture applications appear in Section~\ref{sec:applications}.

\subsection{Related work}
\label{sec:reference-comparison}

On Riemannian three-spheres, Nabutovsky and Rotman bounded the least area
of an embedded minimal two-sphere using sectional curvature, volume, and
diameter \cite{NR2003Sphere}. Their curvature-free bound for stationary
$k$-varifolds gives the recursion and factorial coefficient in Theorem~\ref{thm:mass}
\cite[Theorem~1.4(i)]{NR2006} (see Remark~\ref{mass:NR}).
For closed four-manifolds with vanishing first homology, Wu and Zhu bounded
the first filling function under two-sided Ricci, volume, and diameter bounds
\cite{WZ2020}. Fu and Zhu treated the Einstein case \cite{FZ2026}.
Sections~\ref{sec:closed-fill} and~\ref{sec:relative-wz} build on \cite{WZ2020}.

Struwe constructed a possibly branched free-boundary minimal disk in
domains of $\R^3$ diffeomorphic to a ball \cite{Struwe1984}. Gr\"uter and
Jost obtained an embedded disk in smooth convex bodies in $\R^3$
\cite{GJDisks1986}. Such disks also exist in Riemannian three-balls with
strictly mean-convex boundary
\cite[Theorem~1.2]{SSWZ2026}. Li developed min--max theory for surfaces
in compact three-manifolds with boundary \cite{Li2015}, and De~Lellis
and Ramic treated hypersurfaces with boundary \cite{DLR2018}. Li and Zhou
identified a gap concerning properness in Li's proof \cite[Section~1.1]{LZ2021}.
In ambient dimensions three through seven, they proved existence without
boundary convexity, allowing the embedded hypersurface to touch the boundary
at interior points \cite[Theorem~1.1]{LZ2021}. Further results include
mod-two min--max theory and Morse index bounds \cite{GLWZ2021}, generic
properness \cite{SWZ2024}, and higher-dimensional existence
\cite[Theorem~1.8]{BLPS2026}. Theorem~\ref{thm:mass} gives a mass bound
for stationary varifolds in every codimension.

On a closed $n$-manifold $M$, Nabutovsky and Rotman bounded the least mass
of a nontrivial stationary one-cycle by diameter or filling radius, hence
also by $\Vol(M)^{1/n}$, up to dimensional constants \cite{NR2004}. Rotman
found a stationary geodesic net with one vertex and total length at most
$2\,n!\,\diam(M)$ \cite{Rotman2011}. Orthogonal geodesic chords are stationary
with free boundary. Zhou constructed free-boundary geodesics by min--max
\cite{Zhou2016}. Beach bounded chord lengths on convex Riemannian two-disks
\cite{Beach2026}. For more general
manifolds, see \cite[Theorem~1.1 and Corollaries~4.4--4.5]{BCGRS2025}.
Corollary~4.4 there guarantees orthogonality at one endpoint and, for
convex boundary, at both.

Liokumovich, Nabutovsky, and Rotman swept out a two-sphere by loops based
at any prescribed point, with lengths controlled by area and diameter.
Their bound for contracting a disk boundary also uses its perimeter \cite{LNR2015}.
For a bounded domain $\Omega\subset\R^n$, Guth bounded the $k$-width
by $C(n)\Vol(\Omega)^{k/n}$ \cite{Guth2007}. For a Euclidean ball, the
$k$-width equals the volume of an equatorial $k$-disk \cite{JZ2023}.
The hypersurface width bounds of Glynn-Adey and Liokumovich
\cite{GAL2017} require no low-degree homology vanishing.
Remark~\ref{width:rem:comparison} compares control of one sweepout with
filling estimates for every bounding cycle.

Part~I reviews boundary concentration for free-boundary minimal surfaces
and phase transitions, and the first variation of free-boundary varifolds
\cite[Subsection~1.2]{PartI}.

\subsection{Proof strategy and organization}
\label{sec:proof-overview}

We first describe the proof of Theorem~\ref{thm:mass} in dimension
three. Let $M$ be a three-manifold such that $H_1(M,\partial M;\Zt)$
vanishes, and put $A=\FH_1(2D)$. The theorem gives a nonzero integral
two-varifold $V$ with $\|V\|(M)\le12A$, where $12=3\cdot4$ counts the
edges of a triangle and the faces of a tetrahedron. We suppress small
interpolation and filling errors.

Take a fine triangulation of $(M,\partial M)$. For each interior vertex
$v$, choose a path $\gamma_v$ to a nearest boundary point, of length at
most $D$; put $\gamma_v=0$ at boundary vertices. An edge $e=[u,v]$
determines the chain $z_e=e+\gamma_u+\gamma_v$. Its endpoints cancel
modulo two, except for points of $\partial M$, which vanish relatively.
Thus $z_e$ is a relative one-cycle with $\Mass(z_e)\le2D+o(1)$, even
though the two paths may end at different boundary points, and the
homology assumption provides a relative filling $Q_e$ with
$\dr Q_e=z_e$ and $\Mass(Q_e)\le A+o(1)$. For a triangular face $f$,
correct its geometric chain by its three edge fillings:
\[
 T_f=f+\sum_{e<f}Q_e,\quad
 \dr T_f=0,\quad \Mass(T_f)\le3A+o(1).
\]
The edges cancel, and each vertex path occurs twice. Choices are reused
on shared faces and are zero on boundary simplices.

The next step uses min--max instead of a further homology assumption.
Suppose that no varifold as in Theorem~\ref{thm:mass} has mass below
$12A+\varepsilon$. The contraction alternative of
Lemma~\ref{mass:sublevel} then contracts each $T_f$ to zero through a
path $H_f$ of relative two-cycles of mass at most $3A+o(1)$. For each
tetrahedron $\sigma$, let $c_\sigma$ be the path of small cycles from
$\dr\sigma$ to zero constructed inside $\sigma$ in
Lemma~\ref{mass:simplex-families}, and define
\[
 \Psi_\sigma(t)=c_\sigma(t)+\sum_{f<\sigma}H_f(t).
\]
The four corrected faces sum to $\dr\sigma$, so
$\Psi_\sigma(0)=0=\Psi_\sigma(1)$, and
$\sup_t\Mass(\Psi_\sigma(t))\le4\cdot3A+o(1)$. Summing over all
tetrahedra cancels the paths on shared faces, and the sum represents
$[M,\partial M]\ne0$ under Almgren's isomorphism. At least one loop
$\Psi_\sigma$ is therefore nontrivial. Relative min--max applied to
this single loop (the sum need not satisfy the mass bound) produces the
excluded varifold, a contradiction. Right continuity of $\FH_1$, a
fixed annular selection property, and annular compactness remove
$\varepsilon$ and the errors. In higher dimensions the same face
cancellation produces families with more parameters and the coefficient
$(n+1)!/k!$; this relative version of the construction of
Nabutovsky--Rotman \cite{NR2006} is carried out in
Section~\ref{sec:mass}.

For the relative filling estimates, we copy a bounding relative cycle
into both halves of a controlled smooth double. The seam boundaries
cancel modulo two, and doubling a relative filling shows that the
resulting closed cycle bounds. Cutting a closed filling along the seam
gives two relative fillings; choosing the lighter one supplies the
factor $1/2$, and metric comparison gives the powers of $\Lambda$. A
collar compression justifies the cut even at irregular seam
intersections. In dimension four, body--neck geometry supplies the
closed degree-one and degree-two bounds, while complementary regions
and isoperimetry supply the top-degree bound.

Section~\ref{sec:notation} recalls relative chains and varifolds.
Section~\ref{sec:inputs} contains the min--max theory. Its main
statement, Theorem~\ref{thm:relative-selection}, is proved in
Appendix~\ref{app:selection}. The integrality conclusion uses
\cite[Theorem~1.3]{PartI}. Section~\ref{sec:mass} proves the mass
bound. Sections~\ref{sec:closed-fill} and~\ref{sec:relative-wz}
establish the closed filling estimates in dimension four and transfer
them to manifolds with boundary.
Further applications appear in Section~\ref{sec:applications}.

\subsection{Acknowledgments}

The authors are supported in part by the NSFC Grant 12501065. Z.Z. is
also supported in part by National Key R\&D Program of China
2024YFA1015300.

\subsection{Statements and Declarations}

Data sharing is not applicable to this article as no datasets were
generated or analyzed during the current study. ChatGPT was used to
assist with bibliographic searches, language editing, and initial
drafts of portions of some proofs. Every AI-assisted mathematical
argument was independently reconstructed and checked by the authors,
who reviewed and revised all retained material and take full
responsibility for the manuscript.

\section{Relative chains and varifolds}
\label{sec:notation}

We use integral flat chains with coefficients in $\Zt$ on $(M,\partial M)$ \cite{Federer1969}, equivalently the rectifiable relative-cycle model of Li--Zhou \cite[Section 3]{LZ2021} with mod-two coefficients \cite[Section 3]{GLWZ2021}.
Informally, a relative cycle is a chain whose boundary lies in $\partial M$, and chains supported in $\partial M$ count as zero.
Apart from Definition~\ref{def:noconcentration} and the remark following it, the definitions and conventions of this section are those of \cite[Section~2]{PartI}; we repeat them to keep the paper self-contained.

\begin{definition}
\label{def:relativechains}
Let $\II_j(M;\Zt)$ denote integral $j$-chains supported in $M$, and put
\[
 \IC_j(M,\partial M;\Zt)
 =\II_j(M;\Zt)/\II_j(\partial M;\Zt),
 \quad \dr[T]=[\partial T].
\]
The kernel of $\dr$ is denoted by $\ZC_j(M,\partial M;\Zt)$. For a relative class $T$, define
\[
 \Mass(T)=\inf\{\Mass(\widetilde T):[\widetilde T]=T\}.
\]
For relative cycles, we also use the rectifiable model: finite-mass rectifiable $j$-chains whose boundaries are supported in $\partial M$, identified when their difference is supported in $\partial M$. Each such class has a unique canonical rectifiable representative with zero mass on $\partial M$. Its mass equals the displayed infimum over integral representatives. We omit brackets for relative classes when there is no ambiguity.
\end{definition}

The equivalence and approximation by integral representatives are given in \cite[Lemma 3.1]{GLWZ2021} and \cite[Proposition 2.4]{LMN2018}. A canonical representative may have infinite absolute boundary mass, so the infimum need not be attained. We interpret chain identities using integral representatives before passing to the quotient, and analytic limits using canonical representatives. For Borel $E\subset M$, $\|T\|(E)$ is the canonical representative's mass in $E$; its ordinary varifold is $|T|$.

The relative flat norm is
\begin{equation}\label{eq:flatnorm}
 \Flat(T)=\inf_{T=R+\dr Q}
 \bigl(\Mass(R)+\Mass(Q)\bigr),
\end{equation}
where $R$ and $Q$ are relative chains of dimensions $j$ and $j+1$. The flat and mass distances between two classes are $\Flat(T-S)$ and $\Mass(T-S)$, respectively. In particular,
\begin{equation}\label{eq:flatbasic}
 \Flat(T-S)\le\Mass(T-S),\quad
 \Flat(\dr Q)\le\Mass(Q).
\end{equation}
For a chain $T$ and a measurable set $E$, $T\res E$ denotes restriction. For a Lipschitz map $f$, $f_\#T$ denotes pushforward. The mod-two relative fundamental chain is $\MC$, and $[M,\partial M]\ne0$ denotes its homology class.

\begin{definition}
\label{def:stationarity}
Let $G_k(M)$ be the bundle of unoriented $k$-planes in $TM$. A $k$-varifold is a nonnegative Radon measure on $G_k(M)$; the space of these measures is denoted by $\VC_k(M)$. Its mass measure $\|V\|$ is the pushforward under the bundle projection. An integral varifold is represented by a countably $k$-rectifiable set with positive integer multiplicity almost everywhere.

Let $\mathfrak X_{\tan}(M)$ be the smooth vector fields on $M$ tangent to $\partial M$ along the boundary. The first variation of $V$ is
\begin{equation}\label{eq:firstvariation}
 \delta V(X)=\int_{G_k(M)}\diver_PX(x)\,dV(x,P),
 \quad
 \diver_PX=\sum_{a=1}^{k}g(\nabla_{e_a}X,e_a),
\end{equation}
where $e_1,\ldots,e_k$ is an orthonormal basis of $P$. The varifold is stationary with free boundary if $\delta V(X)=0$ for every $X\in\mathfrak X_{\tan}(M)$. A local statement in an open set $U\subset M$ uses fields compactly supported in $U$.
\end{definition}

For a smooth submanifold $\Sigma$ with outward conormal $\eta_\Sigma$, its first variation is
\begin{equation}\label{eq:smoothvariation}
 \delta|\Sigma|(X)
 =-\int_\Sigma g(\mathbf H_\Sigma,X)\,d\HH^k
  +\int_{\partial\Sigma}g(\eta_\Sigma,X)\,d\HH^{k-1},
\end{equation}
where $\mathbf H_\Sigma$ is the mean curvature vector and $\HH^j$ is $j$-dimensional Riemannian Hausdorff measure. A smooth free-boundary minimal submanifold has $\mathbf H_\Sigma=0$ and meets $\partial M$ orthogonally along its boundary.

We use the free-boundary version of Pitts' almost minimizing condition \cite{Pitts1981} developed in \cite[Definitions~3.18 and~4.19]{LZ2021}. For the mod-two hypersurface case, see \cite[Section~4]{GLWZ2021}.

\begin{definition}
\label{def:almostminimizing}
Fix a smooth Riemannian metric on the compact Grassmann bundle $G_k(M)$, and let $d_G$ be its distance. For $V,W\in\VC_k(M)$, define
\begin{equation}\label{eq:varifold-distance}
 \dV(V,W)=\sup_{\substack{\|f\|_\infty\le1\\\Lip_{d_G}(f)\le1}}
 \left|\int_{G_k(M)}f\,dV-\int_{G_k(M)}f\,dW\right|.
\end{equation}
The supremum is over real-valued functions on $G_k(M)$ whose supremum norm and Lipschitz constant with respect to $d_G$ are at most one. This distance induces weak convergence on sets of varifolds with bounded mass.
Given an open set $U\subset M$, $\eps,\delta>0$, and $d\in\{\Flat,\Mass\}$, let $\AMset(U;\eps,\delta;d)$ be the relative $k$-cycles $T$ such that every finite sequence $T=T_0,\ldots,T_N$ satisfying
\[
 \spt(T_i-T)\subset U,\quad
 d(T_i-T_{i-1})\le\delta,\quad
 \Mass(T_i)\le\Mass(T)+\delta
\]
also satisfies $\Mass(T_N)\ge\Mass(T)-\eps$. Supports are taken for canonical representatives.

A varifold $V$ is almost minimizing modulo two in $U$ if, for every $\eps>0$, there are $\delta>0$ and $T\in\AMset(U;\eps,\delta;\Flat)$ such that $\dV(|T|,V)<\eps$.

Choose a smooth extension of the metric across $\partial M$. For balls in this extension, define the relative annulus
\[
 \An(p;s,r)=M\cap\bigl(B_r(p)\setminus\overline B_s(p)\bigr),
 \quad0<s<r.
\]
The varifold is almost minimizing in small relative annuli if, for each $p\in M$, there is $r_p>0$ such that it is almost minimizing in every $\An(p;s,r)$ with $0<s<r<r_p$.
\end{definition}

In words, a cycle $T$ belongs to $\AMset(U;\eps,\delta;d)$ if its mass cannot be decreased by more than $\eps$ through a finite chain of modifications inside $U$, each of which is $\delta$-small and none of which raises the mass by more than $\delta$. A varifold is almost minimizing if it is approximated by such cycles with $\eps\to0$. The cycles and the steps are relative, whereas $|T|$ and $V$ are ordinary varifolds on $M$.

Canonical representatives have zero boundary mass, although their ordinary varifold limits may acquire it. For parameter families we also use the following condition.

\begin{definition}
\label{def:noconcentration}
Let $X$ be a finite cubical complex. A flat-continuous map $\Phi:X\to\ZC_k(M,\partial M;\Zt)$ has no concentration of mass if
\begin{equation}\label{eq:noconcentration}
 \lim_{r\downarrow0}\sup_{x\in X,\,p\in M}
 \|\Phi(x)\|(B_r(p))=0.
\end{equation}
\end{definition}

Every mass-continuous map from a compact parameter space has no
concentration of mass for $k\ge1$; this is
\cite[Lemma 2.7]{LMN2018}, applied to the canonical representatives
of the relative cycles.

We finally record the compactness and small-filling properties of relative cycles.

\begin{lemma}
\label{inp:compactness}
Let $1\le k\le n-1$. Relative cycles with a common mass bound are sequentially compact in the relative flat topology, and relative mass is lower semicontinuous. There are constants $\rho_M,C_M>0$ such that, if
$\Flat(T-S)<\rho_M$, then there is a relative $(k+1)$-chain $Q$ satisfying
\begin{equation}\label{inp:eq:small-filling}
 \dr Q=T-S,
 \quad\Mass(Q)\le C_M\Flat(T-S).
\end{equation}
\end{lemma}

Compactness and lower semicontinuity follow from
\cite[Section~2.3 and Proposition~2.4]{LMN2018}. In the compactness
argument, projection in a boundary collar gives equivalent integral
representatives with uniformly bounded boundary mass, to which
Federer--Fleming compactness applies.

For the small-filling estimate, apply Almgren's relative isoperimetric
choice \cite[Proposition~1.23]{Almgren1962} with mod-two coefficients,
as in \cite[Appendix~A, p.~957]{LMN2018}. For integral representatives
$\widetilde T,\widetilde S$ of the two cycles, it gives a mod-two chain
$\widetilde Q$ in $M$ of mass at most $C_M\Flat(T-S)$ such that
$\partial\widetilde Q-(\widetilde T-\widetilde S)$ is supported in
$\partial M$. Its relative class $Q=[\widetilde Q]$ therefore satisfies
\eqref{inp:eq:small-filling}.

\section{Relative interpolation and min--max theory}
\label{sec:inputs}

This section collects the tools used in the filling construction. Theorem~\ref{inp:interpolation} provides interpolation for relative cycles; Lemma~\ref{inp:localization} compares the mass and flat versions of almost minimization; and Theorem~\ref{inp:almgren} is Almgren's isomorphism. From these we derive two statements that are used repeatedly. First, a family of cycles of small mass can be contracted (Lemma~\ref{mass:small-cone}). Second, a family of cycles either can be contracted below a prescribed mass ceiling or produces a stationary varifold below that ceiling (Theorem~\ref{thm:relative-selection}); the second alternative is proved in Appendix~\ref{app:selection}. Finally, Lemma~\ref{inp:annular-compactness} allows us to pass to limits of such varifolds. The dimensions of the parameter space and of the cycles are independent. We assume $2\le k\le n-1$ throughout, and no hypersurface regularity is used.

Write $I=[0,1]$. If $X$ is a finite cubical complex, write $X(j)_0$ for the vertices of its subdivision into cubes of side length $3^{-j}$. For a vertex map $\phi$, define its mass fineness by
\[
 \mathbf f(\phi)=\max\{\Mass(\phi(v)-\phi(w)):
       v,w\text{ are endpoints of one subdivided edge}\}.
\]
We set $\mathbf f(\phi)=0$ if the subdivision has no edges.
The following results of Liokumovich--Marques--Neves concern general-dimensional mod-two relative cycles.

\begin{theorem}
\label{inp:interpolation}
The following statements hold for every fixed finite parameter dimension.
\begin{enumerate}[label=\textup{(\roman*)},leftmargin=2.2em]
\item If $\Phi:X\to\ZC_k(M,\partial M;\Zt)$ is flat-continuous, has bounded mass, and has no concentration of mass, there are subdivisions $j_i\to\infty$, numbers $e_i\downarrow0$, and vertex maps $\phi_i:X(j_i)_0\to\ZC_k(M,\partial M;\Zt)$ such that
\[
 \mathbf f(\phi_i)<e_i,\quad
 \sup_v\Flat(\phi_i(v)-\Phi(v))<e_i,\quad
 \max_v\Mass(\phi_i(v))\le\sup_x\Mass(\Phi(x))+e_i.
\]
If $\Phi$ is mass-continuous, one may take
$\phi_i$ to be its restriction to the vertices of a sufficiently fine
subdivision. This retains prescribed vertex values and zero parameter faces.
\item There are constants $\delta_0,C_0>0$, independent of $j$, such that every map
\[
 \phi:X(j)_0\to\ZC_k(M,\partial M;\Zt),\quad \mathbf f(\phi)<\delta_0,
\]
admits a mass-continuous extension $\Psi:X\to\ZC_k(M,\partial M;\Zt)$ satisfying
\begin{equation}\label{inp:eq:almgren-extension}
 \Psi(v)=\phi(v)\quad(v\in X(j)_0),\quad
 \sup_{x,y\in\alpha}\Mass\bigl(\Psi(x)-\Psi(y)\bigr)
 \le C_0\mathbf f(\phi)
\end{equation}
for every cell $\alpha$ of $X(j)$. The restriction of $\Psi$ to $\alpha$ depends only on the values of $\phi$ at the vertices of $\alpha$. This map is called the Almgren extension of $\phi$. For the maps $\phi_i$ in \textup{(i)}, their Almgren extensions are homotopic to the original map $\Phi$ in \textup{(i)} through flat-continuous maps for all sufficiently large $i$.
\end{enumerate}
\end{theorem}

See \cite[Theorems 2.10 and 2.11, Proposition 2.12, and Appendices A and B]{LMN2018}.
The compactness and interpolation proofs there do not use the
higher-homology assumption imposed in their Section~2.1 for the subsequent
spectrum construction.

\begin{lemma}
\label{inp:localization}
Let $K\subset U\subset M$, where $K$ is compact and $U$ is relatively open, and fix $\tau\in\ZC_k(M,\partial M;\Zt)$. Given $L,\eta,\beta>0$, there is $\epsilon_0>0$ such that any two relative cycles $S_1,S_2$ satisfying
\[
 \Mass(S_i)\le L,\quad
 \operatorname{spt}(S_i-\tau)\subset K,\quad
 \Flat(S_1-S_2)<\epsilon_0,
\]
can be joined by a finite sequence $S_1=R_0,R_1,\ldots,R_N=S_2$ with
\[
 \begin{gathered}
 \spt(R_a-\tau)\subset U,\quad \Mass(R_a)\le L+\eta
 \quad(0\le a\le N),\\
 \Mass(R_a-R_{a-1})\le\beta\quad(1\le a\le N).
 \end{gathered}
\]
The number $\epsilon_0$ may depend on $L,\eta,\beta,K,U$, and $\tau$. In particular, if $W\Subset U$ and $T\in\AMset(U;\zeta,\delta;\Mass)$, then for some $\delta'>0$,
\[
 T\in\AMset(W;\zeta,\delta';\Flat).
\]
Thus a varifold that is almost minimizing using mass steps in $U$ is almost minimizing using flat steps in $W$.
\end{lemma}

See \cite[Lemma B.1 and the proof of Theorem 3.20 in Appendix B]{LZ2021}, whose restriction and filling arguments also apply with mod-two coefficients, and \cite[Lemma 3.8 and Theorem 3.9]{Pitts1981} for the interior case. To separate the mass allowance $\eta$ from the step bound $\beta$, use that lemma with its parameter smaller than $\min\{\eta,\beta\}$.

For the stated consequence, apply the first part with $K=\overline W$, $\tau=T$, $L=\Mass(T)+\delta/2$ and $\eta=\beta=\delta/2$, and choose $\delta'<\min\{\delta/2,\epsilon_0\}$. Let $T=T_0,\ldots,T_N$ be a sequence with $\spt(T_i-T)\subset W$, flat steps at most $\delta'$ and masses at most $\Mass(T)+\delta'$. Interpolating between consecutive terms gives a sequence from $T$ to $T_N$ whose differences from $T$ are supported in $U$, with mass steps at most $\delta$ and masses at most $\Mass(T)+\delta$. Hence $\Mass(T_N)\ge\Mass(T)-\zeta$.

\begin{theorem}
\label{inp:almgren}
For $m\ge1$, there is an isomorphism
\begin{equation}\label{inp:eq:almgren}
 \Alm:\pi_m\bigl(\ZC_k(M,\partial M;\Zt),0\bigr)
 \longrightarrow H_{k+m}(M,\partial M;\Zt),
\end{equation}
where the cycle space has the relative flat topology. The corresponding description by sufficiently fine discrete families gives the same class. The map is additive and is compatible with taking boundaries of parameter cells.
\end{theorem}

For Almgren's isomorphism theorem and its relative mod-two formulation, see \cite{Almgren1962} and \cite[Section 2]{LMN2018}, respectively.

For $r\ge1$, a based $r$-sphere of cycles will be represented by a map on $I^r$ that is zero on $\partial I^r$. For $r=0$, $I^0$ is a single point and $\partial I^0$ is empty.

\begin{lemma}\label{mass:small-cone}
Let $X$ be a finite cubical complex and $A\subset X$ a subcomplex. There are constants $a_X,C_X>0$, depending only on $M,g,k$ and $(X,A)$, such that every mass-continuous map $f:X\to\ZC_k(M,\partial M;\Zt)$ with $f|_{A}=0$ and $L:=\sup_X\Mass(f)<a_X$ admits a mass-continuous homotopy $H:X\times I\to\ZC_k(M,\partial M;\Zt)$ satisfying
\[
 H(\cdot,0)=f,\quad H(\cdot,1)=0,\quad
 H|_{A\times I}=0,\quad \sup_{X\times I}\Mass(H)\le C_XL.
\]
For $(X,A)=(I^r,\partial I^r)$, this is an extension to $I^{r+1}$ that is zero on every boundary face except the bottom.

The same conclusion holds for a family $f$ of absolute cycles supported in a convex simplex $S$ in a fixed Euclidean coordinate ball. Here continuity is in the absolute mass topology, $L=\sup_{x\in X}\Mass_{\mathrm{abs}}(f(x))$, and the homotopy consists of absolute cycles supported in $S$ with absolute mass at most $C_XL$. Its constants may depend on the coordinate ball but are independent of the size and shape of $S$.
\end{lemma}
\begin{proof}
The idea is to subtract from $f$ a discretized copy of itself, which can be contracted by an Almgren extension, and to iterate: the error left after each step has less than half the mass of the previous one.

We first describe one step of this iteration.
Throughout we use the cellwise mass extension of Theorem~\ref{inp:interpolation}\textup{(ii)} and \cite[Theorem~2.11]{LMN2018}; we make the same choices on shared faces, and we choose the zero extension on every subcomplex on which the vertex data vanish. Let $\delta_0,C_0$ be extension constants valid for both $X$ and $X\times I$, and choose $a_X<\delta_0/2$.

Let $g:X\to\ZC_k(M,\partial M;\Zt)$ be mass-continuous with $g|_A=0$ and $0<a:=\sup_X\Mass(g)\le L$. Sample $g$ at the vertices of a fine subdivision of $X$, and let $F$ be the Almgren extension of these samples. If the mass oscillation of $g$ on every cell of the subdivision is at most $\omega$, then comparing $g(x)$ and $F(x)$ with their common value at a vertex of a cell containing $x$ gives
$\sup_X\Mass(g-F)\le(1+C_0)\omega$.
Since $g$ is uniformly mass-continuous, a sufficiently fine subdivision makes this less than $a/2$.

Next, on the corresponding subdivision of $X\times I$, assign the samples of $g$ to the vertices of the bottom $X\times\{0\}$ and zero to the vertices of every higher layer. The resulting vertex map $\phi$ satisfies $\max_v\Mass(\phi(v))\le a$ and $\mathbf f(\phi)\le2a<\delta_0$. Extend it over the product, using on the bottom the cellwise choices already made for $F$. The resulting Almgren extension $J$ satisfies $J(\cdot,0)=F$ and
\[
 \sup\Mass(J)\le a+C_0\mathbf f(\phi)
 \le Ca,\quad C=1+2C_0,
\]
and it vanishes on $X\times\{1\}$ and on $A\times I$, where all vertex data vanish. Put $g'=g-F$. Then $(x,t)\mapsto g'(x)+J(x,t)$ is a mass-continuous homotopy which starts at $g'+F=g$, ends at $g'$, vanishes on $A\times I$, and has mass at most $(C+1)a$; moreover $\sup_X\Mass(g')<a/2$.

We now iterate this construction.
Set $f_0=f$ and $a_\ell=\sup_X\Mass(f_\ell)$. If $a_\ell=0$, stop. Otherwise apply the construction above to $g=f_\ell$ and put $f_{\ell+1}=g'$. Then $a_\ell\le2^{-\ell}L$, and the $\ell$th homotopy joins $f_\ell$ to $f_{\ell+1}$ with mass at most $(C+1)2^{-\ell}L$. Run the $\ell$th homotopy on the time interval $[1-2^{-\ell},1-2^{-\ell-1}]$ and set $H(\cdot,1)=0$; if the iteration stops at $f_\ell=0$, let $H$ vanish after that time. Adjacent pieces have the same endpoint map, and the bound $(C+1)2^{-\ell}L$, which is uniform in $x$, gives joint mass continuity at time one. Thus $H$ is a mass-continuous homotopy from $f$ to zero which vanishes on $A\times I$ and satisfies $\sup\Mass(H)\le C_XL$ with $C_X=C+1$.

For the last assertion, let $\overline B$ be the fixed closed coordinate ball containing $S$. Embed a neighborhood $B'$ of $\overline B$ into a round sphere $N$ by a map $\iota$, and extend the inverse coordinates to a smooth map $q:N\to\R^n$ using a cutoff equal to one near $\iota(\overline B)$. Let $\pi_S$ be the nearest-point projection onto the convex set $S$, and set $P=\pi_S\circ q$. Then $P\circ\iota=\mathrm{id}$ on $S$ and $\Lip(P)\le\Lip(q)$, because $\pi_S$ is $1$-Lipschitz. Apply the construction above in $N$, using the absolute-cycle case of \cite[Theorem 2.11]{LMN2018}, to obtain a contraction $\widehat H$ of $\iota_\#f$, and put $H=P_\#\widehat H$. Pushforward by $P$ preserves cycles and zero values and is Lipschitz for the mass distance; since $f$ is supported in $S$, it returns the initial map $P_\#\iota_\#f=f$ exactly. Thus $H$ is a contraction of $f$ through cycles supported in $S$, and
$\sup\Mass(H)\le C_N\Lip(q)^k\Lip(\iota|_B)^kL$.
Both the coefficient and the smallness threshold depend on the coordinate ball but not on $S$.
\end{proof}

The min--max theorem used in Section~\ref{sec:mass} keeps the homotopies
below a mass ceiling $\Lambda$. The resulting varifold satisfies the
following annular property, which persists under varifold convergence.

Fix $R_*>0$ smaller than one quarter of the common coordinate radius used for local interpolation. For a positive integer $N$, write $\mathcal P_N$ for the following property:
\begin{equation}\label{inp:eq:annular-property}
\begin{gathered}
 \text{for every }p\in M\text{ and every collection }
 U_j=\An(p;s_j,r_j),\quad 1\le j\le N,\\
 0<s_j<r_j<R_*,\quad s_j>12r_{j+1}\quad(1\le j<N),\\
 V\text{ is almost minimizing in at least one }U_j.
\end{gathered}
\end{equation}
In words, $V$ cannot fail to be almost minimizing in each of $N$ well-separated concentric relative annuli of outer radius less than $R_*$. The radius $R_*$ remains fixed throughout. The same property, at a fixed scale, is the hypothesis of the depth estimate \cite[Corollary~1.8]{PartI}. For the related closed-hypersurface condition and its use under limits, see \cite[Definition 3.3 and proof of Lemma 3.20]{CL2025}. We give the relative-cycle argument.

The annular property implies the two properties of the varifold that are used later.

\begin{lemma}\label{inp:PN-consequences}
Let $V$ be a $k$-varifold in $M$ that is stationary with free boundary and satisfies $\mathcal P_N$ for some positive integer $N$. Then $V$ is almost minimizing in small relative annuli, and $V$ is integral.
\end{lemma}

\begin{proof}
Fix $p\in M$, and suppose that for every $\rho>0$ there is an annulus $\An(p;s,r)$ with $0<s<r<\rho$ in which $V$ is not almost minimizing. Then we can choose successively annuli $U_j=\An(p;s_j,r_j)$, $1\le j\le N$, with $r_1<R_*$ and $r_{j+1}<s_j/12$, in none of which $V$ is almost minimizing. This contradicts $\mathcal P_N$. Hence $V$ is almost minimizing in small relative annuli, and it is integral by \cite[Corollary~3.5]{PartI}.
\end{proof}

\begin{theorem}\label{thm:relative-selection}
Let $X$ be a cubical subcomplex of a finite subdivision of $I^m$, $m\ge0$, and let $A\subset X$ be a subcomplex. Suppose $\Phi_0:X\to\ZC_k(M,\partial M;\Zt)$ is mass-continuous and zero on $A$. Fix
$\Lambda>\sup_{x\in X}\Mass(\Phi_0(x))$.
Let $\mathcal C$ consist of the maps joined to $\Phi_0$, relative to $A$, by mass-continuous homotopies whose values have mass less than $\Lambda$, and put
\begin{equation}\label{inp:eq:restricted-width}
 a=\inf_{\Phi\in\mathcal C}\sup_{x\in X}\Mass(\Phi(x)).
\end{equation}
There is a positive integer $N(m)$, depending only on $m$, for which:
\begin{enumerate}[label=\textup{(\roman*)},leftmargin=2.2em]
\item If $a=0$, there is a mass-continuous homotopy $H:X\times I\to\ZC_k(M,\partial M;\Zt)$ such that
\[
 H(\cdot,0)=\Phi_0,\quad H(\cdot,1)=0,\quad
 H|_{A\times I}=0,\quad \sup_{X\times I}\Mass(H)<\Lambda.
\]
\item If $a>0$, there is an integral $k$-varifold $V$, stationary with free boundary and satisfying $\mathcal P_{N(m)}$, such that
\[
 0<\|V\|(M)=a\le\sup_{x\in X}\Mass(\Phi_0(x))<\Lambda.
\]
In particular, $V$ is almost minimizing in small relative annuli.
\end{enumerate}
One may take $N(m)=2\cdot3^m-1$, which is nondecreasing in $m$. The case $m=0$ includes a single cycle, with $X$ a point and $A=\varnothing$.
\end{theorem}

\begin{proof}
If $a=0$, take $\Phi\in\mathcal C$ with
$\sup\Mass(\Phi)<\min\{a_X,\Lambda/(2C_X)\}$, where $a_X$ and $C_X$ are the constants of Lemma~\ref{mass:small-cone}. That lemma contracts $\Phi$ to zero, relative to $A$, through maps of mass less than $\Lambda/2$. Concatenating this contraction with a homotopy from $\Phi_0$ to $\Phi$ as in the definition of $\mathcal C$ gives \textup{(i)}.

For $a>0$ we adapt the argument of Pitts \cite[Theorem~4.10]{Pitts1981}
and its free-boundary version for integer hypersurface cycles due to
Li and Zhou \cite[Theorem~4.21]{LZ2021}. Appendix~\ref{app:selection}
carries out the selection for relative mod-two $k$-cycles in every
dimension $2\le k\le n-1$, keeping the homotopies below $\Lambda$ and
obtaining the fixed-scale property $\mathcal P_{N(m)}$.
A sequence $\Phi_i\in\mathcal C$ with $\sup_X\Mass(\Phi_i)\to a$ is pulled tight, so that the varifold limits of its cycles of nearly maximal mass are stationary with free boundary. If none of these limits satisfied $\mathcal P_{N(m)}$, each such cycle could be deformed inside any one of $N(m)$ concentric annuli so that its mass decreases by a definite amount (Lemma~\ref{inp:localization}). The number $N(m)$ is large enough to select these annuli so that neighboring parameter vertices use disjoint ones (Lemma~\ref{app:selection-lemma}). The deformations can then be performed simultaneously (Lemma~\ref{app:combination}), and they yield an element of $\mathcal C$ of maximal mass less than $a$, which is impossible. A limit $V$ satisfying $\mathcal P_{N(m)}$ has mass $a$, and it is integral and almost minimizing in small relative annuli by Lemma~\ref{inp:PN-consequences}.
\end{proof}

For based families this gives the following alternative.

\begin{lemma}\label{mass:sublevel}
Let $r\ge0$, and let $f:I^r\to\ZC_k(M,\partial M;\Zt)$ be mass-continuous and zero on $\partial I^r$. Write $L=\sup\Mass(f)$. For every $e>0$, at least one of the following conclusions holds:
\begin{enumerate}[label=\textup{(\roman*)},leftmargin=2.2em]
\item There is a nonzero integral $k$-varifold, stationary with free boundary and almost minimizing modulo two in small relative annuli, whose mass is at most $L$. It can be chosen to satisfy $\mathcal P_{N(r)}$, with $N(r)$ as in Theorem~\ref{thm:relative-selection}.
\item There is a mass-continuous extension $H:I^{r+1}\to\ZC_k(M,\partial M;\Zt)$ equal to $f$ on the bottom face and zero on all other boundary faces, with supremum of mass less than $L+e$.
\end{enumerate}
For $r=0$, the second conclusion is a path from the single cycle to zero.
\end{lemma}
\begin{proof}
Apply Theorem~\ref{thm:relative-selection} with
$(X,A)=(I^r,\partial I^r)$, $\Phi_0=f$, and $\Lambda=L+e$.
For $r=0$, take $X$ to be a point and $A=\varnothing$.
\end{proof}

The positive-width alternative applies in particular to a family with nonzero Almgren class.

\begin{corollary}\label{inp:sweepout}
Let $m=n-k$ and let
$\Phi:I^m\to\ZC_k(M,\partial M;\Zt)$ be mass-continuous, with
$\Phi|_{\partial I^m}=0$ and
\[
 \Alm([\Phi])\ne0\quad\text{in }H_n(M,\partial M;\Zt).
\]
Then there is a nonzero integral $k$-varifold $V$, stationary with free boundary and almost minimizing in small relative annuli, such that
\[
 0<\|V\|(M)\le\sup_{x\in I^m}\Mass(\Phi(x)).
\]
Moreover, $V$ can be chosen to satisfy $\mathcal P_{N(m)}$.
\end{corollary}
\begin{proof}
Apply Theorem~\ref{thm:relative-selection} to
$(X,A)=(I^m,\partial I^m)$, with $\Phi_0=\Phi$ and any
$\Lambda>\sup\Mass(\Phi)$.
Alternative \textup{(i)} would give a based mass-continuous, hence flat-continuous, null-homotopy of $\Phi$, contradicting
$\Alm([\Phi])\ne0$ by Theorem~\ref{inp:almgren}.
Alternative \textup{(ii)} gives the conclusion.
\end{proof}

\begin{lemma}
\label{inp:annular-compactness}
Fix a positive integer $N$. Let $V_i$ be nonzero integral $k$-varifolds, stationary with free boundary, with uniformly bounded masses and satisfying $\mathcal P_N$. A subsequence converges to a nonzero integral varifold $V$, stationary with free boundary and satisfying $\mathcal P_N$, with
\[
 \|V\|(M)=\lim_i\|V_i\|(M).
\]
In particular, $V$ is almost minimizing in small relative annuli.
\end{lemma}
\begin{proof}
After passing to a subsequence, $V_i\rightharpoonup V$. Total masses converge since $M$ is compact, and stationarity passes to the limit in \eqref{eq:firstvariation}. It remains to show that $V$ satisfies $\mathcal P_N$, that it is integral and almost minimizing in small annuli, and that it is nonzero.

We begin with the annular property. The key observation is that almost minimization in a fixed open set $U$ is closed under varifold convergence: if $W_i\rightharpoonup W$ and every $W_i$ is almost minimizing in $U$, then for any $\varepsilon>0$ choose $i,T,\delta$ with
\[
 \dV(W_i,W)<\varepsilon/2,\quad
 \dV(|T|,W_i)<\varepsilon/2,\quad
 T\in\AMset(U;\varepsilon/2,\delta;\Flat).
\]
Then $\dV(|T|,W)<\varepsilon$ and $T\in\AMset(U;\varepsilon,\delta;\Flat)$, so $W$ is almost minimizing in $U$. Now consider any collection of annuli as in \eqref{inp:eq:annular-property}. Each $V_i$ is almost minimizing in at least one of its $N$ annuli, so one fixed annulus has this property for infinitely many $i$, and hence for $V$. The set of indices may depend on the chosen collection, but the limit is the already fixed $V$, so no simultaneous diagonal choice over centers and radii is required. Thus $V$ satisfies $\mathcal P_N$. Almost minimization in small annuli and integrality now follow from Lemma~\ref{inp:PN-consequences}.

It remains to prove that $V$ is nonzero. For this we show that every nonzero integral varifold $W$ which is stationary with free boundary satisfies $\|W\|(M)\ge c'\rho^k$, where $c',\rho>0$ depend only on $M,g,k$. Applied to the $V_i$, this gives $\|V\|(M)=\lim_i\|V_i\|(M)>0$. Fix an isometric embedding $M\subset\R^q$ and a uniform radius $\rho>0$ for interior and boundary monotonicity \cite[Section 5]{Allard1972}, \cite[Theorem 2.3]{LZ2021}. Choose a point $x$ at which $W$ has positive integer density and put $d=\dist_{\R^q}(x,\partial M)$. If $d=0$ or $d\ge\rho/4$, the corresponding monotonicity estimate gives $\|W\|(M)\ge c\rho^k$. Otherwise, let $p$ be a nearest boundary point and write $B_r^e$ for Euclidean balls. Interior monotonicity at $x$ bounds the mass of $B_{d/2}^e(x)\subset B_{2d}^e(p)$ from below, and boundary monotonicity at $p$ then gives
\[
 \|W\|(M)\ge c\left(\frac{\rho}{2d}\right)^k
 \|W\|(B_{2d}^e(p))
 \ge c^2\left(\frac{\rho}{2d}\right)^k(d/2)^k
 =c^2\,4^{-k}\rho^k.
\]
Here $c>0$ depends only on $M,g,k$.
\end{proof}

\section{Mass estimates via relative filling functions}
\label{sec:mass}

Throughout this section we assume \eqref{eq:homology}, recall $D$ and $b_j$ from \eqref{eq:boundarydistance} and~\eqref{eq:budgets}, and set $A=b_{k-1}$. The proof of Theorem~\ref{thm:mass} has two stages. For every $\eta>0$, we first construct a nonzero integral $k$-varifold, stationary with free boundary and almost minimizing modulo two in small relative annuli, of mass less than
\begin{equation}\label{mass:limiting-bound}
 \frac{(n+1)!}{k!}A+\eta.
\end{equation}
All these varifolds satisfy one and the same annular property $\mathcal P_N$, so Lemma~\ref{inp:annular-compactness} then gives a limit satisfying the exact bound.

The construction takes place on a fine triangulation of $(M,\partial M)$ and proceeds by induction on the dimension of the simplices. Below dimension $k$, the homology assumption allows us to fill cycles (Lemma~\ref{mass:chain-homotopy}); this is where the filling functions enter. From dimension $k$ on, no homology assumption is available, and we argue by contradiction instead: if there were no varifold of small mass, then Lemma~\ref{mass:sublevel} would contract every family of cycles attached to a simplex, and these contractions could be assembled, simplex by simplex, up to dimension $n$. The obstruction to completing this induction is the fundamental class. Indeed, Theorem~\ref{inp:almgren}, with the flat topology, gives
\begin{equation}\label{mass:almgren-source}
 \pi_{n-k}\bigl(\ZC_k(M,\partial M;\Zt),0\bigr)
 \cong H_n(M,\partial M;\Zt)
 \ni[M,\partial M]\ne0,
\end{equation}
and the families attached to the $n$-simplices add up to a family representing $[M,\partial M]$. One of them is therefore homotopically nontrivial, and Corollary~\ref{inp:sweepout} produces the varifold after all. No assumption on $H_k(M,\partial M;\Zt)$ or $H_{k+1}(M,\partial M;\Zt)$ is needed.

We first establish the properties of the filling functions that are needed to remove the errors.

\begin{lemma}\label{mass:filling-finite}
Suppose $1\le j<n$ and $H_j(M,\partial M;\Zt)=0$. Then $\FH_j(0)=0$, and $0<\FH_j(v)<\infty$ for every $v>0$. The function $\FH_j$ is nondecreasing and right-continuous on $[0,\infty)$.
\end{lemma}
\begin{proof}
Write $\operatorname{Fill}_j(z)$ for the infimum of the masses of relative
fillings of $z$. The relative deformation theorem identifies current
homology with singular relative homology, so the homology hypothesis
makes this infimum finite; see \cite{Federer1969,Almgren1962}.
Adding a small filling of $z'-z$ to a filling of $z$, and reversing the
roles of $z,z'$, gives from Lemma~\ref{inp:compactness}
\begin{equation}\label{mass:fill-continuity}
 |\operatorname{Fill}_j(z')-\operatorname{Fill}_j(z)|
 \le C_j\Flat(z'-z)
 \quad\text{if }\Flat(z'-z)<\rho_j,
\end{equation}
where $C_j,\rho_j>0$ are its constants in degree $j$.
Thus $\operatorname{Fill}_j$ is continuous on the compact flat space of
cycles of mass at most $v$, and $\FH_j(v)$ is finite.

Monotonicity is immediate. If $v_\ell\downarrow v$, choose cycles with
\begin{equation}\label{mass:near-supremum}
 \Mass(z_\ell)\le v_\ell,\quad
 \FH_j(v_\ell)-\ell^{-1}\le\operatorname{Fill}_j(z_\ell)
 \le\FH_j(v_\ell).
\end{equation}
A subsequence converges flatly to $z$ with $\Mass(z)\le v$.
Equation~\eqref{mass:fill-continuity} then gives
$\lim_\ell\FH_j(v_\ell)=\operatorname{Fill}_j(z)\le\FH_j(v)$;
the reverse inequality follows from monotonicity. This proves right
continuity, including at zero. Finally, $\FH_j(0)=0$, and for $v>0$ the
boundary of a sufficiently small interior $(j+1)$-disk is a nonzero
relative cycle $z$ of mass at most $v$. Since every filling has mass at
least $\Flat(z)>0$, we have $\FH_j(v)>0$.
\end{proof}

By Lemma~\ref{mass:filling-finite}, each $b_j$ is finite and positive.
Fix $s>0$. We subdivide the triangulation until each simplex of positive
dimension has mass less than $s$, and choose each filling with mass
less than its infimal mass plus $s$. With these allowances, set
\begin{equation}\label{mass:padded-fillings}
 \widehat b_1(s)=\FH_1(2D+s)+s,\quad
 \widehat b_j(s)=\FH_j((j+1)\widehat b_{j-1}(s)+s)+s
 \quad(2\le j\le k-1).
\end{equation}
Set $A_s=\widehat b_{k-1}(s)$. By monotonicity and right continuity, $\widehat b_j(s)\downarrow b_j$ inductively as $s\downarrow0$, so $A_s\to A$.

For the families of cycles attached to simplices of dimension $j\ge k$, the mass bounds will be
\begin{equation}\label{mass:parameter-budgets}
 L_k(s)=(k+1)A_s+2s,\quad
 L_j(s)=(j+1)L_{j-1}(s)+2s\quad(k+1\le j\le n).
\end{equation}
The two errors account for the geometric family and Lemma~\ref{mass:sublevel}. Iteration gives
\begin{equation}\label{mass:budget-formula}
 L_j(s)=\frac{(j+1)!}{k!}A_s
       +2s\sum_{\ell=k}^{j}\frac{(j+1)!}{(\ell+1)!}
 \quad(k\le j\le n).
\end{equation}
In particular,
\begin{equation}\label{mass:budget-limit}
 \lim_{s\downarrow0}L_n(s)=\frac{(n+1)!}{k!}A.
\end{equation}
For each $s>0$ the construction is finite; passage to a varifold limit will use Lemma~\ref{inp:annular-compactness}.

Fix $s>0$ and a fine smooth triangulation $K$ of $(M,\partial M)$ \cite{Whitney1957}. Write $K_j$ for its $j$-simplices and $\tau<\sigma$ for the face relation. Let $\widehat\iota_j(\sigma)=\llbracket\sigma\rrbracket$ and let $\iota_j(\sigma)$ be its relative class, zero for boundary simplices. Before choosing fillings, refine so that every positive-dimensional simplex has mass less than $s$ and Lemma~\ref{mass:simplex-families} applies.

\begin{lemma}\label{mass:chain-homotopy}
For every $\sigma\in K_j$, $0\le j\le k-1$, one can choose a relative $(j+1)$-chain $h_j(\sigma)$ with the following properties:
\begin{align}
 \dr h_0(v)&=\iota_0(v),&\Mass(h_0(v))&\le D,
 \label{mass:vertex}\\
 \dr h_j(\sigma)&=\iota_j(\sigma)
       +\sum_{\substack{\tau<\sigma\\\dim\tau=j-1}}h_{j-1}(\tau),
 &\Mass(h_j(\sigma))&<\widehat b_j(s)
 \quad(1\le j\le k-1).
 \label{mass:chain-identity}
\end{align}
All these chains are zero on boundary simplices. For each $k$-simplex $f$, the chain
\begin{equation}\label{mass:initial-cycle}
 T_f=\iota_k(f)+
       \sum_{\substack{\tau<f\\\dim\tau=k-1}}h_{k-1}(\tau)
\end{equation}
is a relative $k$-cycle and satisfies
\begin{equation}\label{mass:initial-budget}
 \Mass(T_f)<(k+1)A_s+s=L_k(s)-s.
\end{equation}
\end{lemma}
\begin{proof}
For each interior vertex $v$, let $h_0(v)$ be a shortest path to
$\partial M$, regarded as a relative chain, and put $h_0(v)=0$ on boundary
vertices. Then $\dr h_0(v)=\iota_0(v)$ and $\Mass(h_0(v))\le D$.
The paths may end at different boundary points.

Inductively, for a $j$-simplex $\sigma$, set
\[
 z_\sigma=\iota_j(\sigma)+
 \sum_{\substack{\tau<\sigma\\\dim\tau=j-1}}h_{j-1}(\tau).
\]
This is a relative cycle: the geometric boundary terms cancel, and
every codimension-two face occurs twice. For $j=1$ its mass is less
than $2D+s$; for $j\ge2$ it is less than
$(j+1)\widehat b_{j-1}(s)+s$. Choose a filling within $s$ of its infimal
mass. By \eqref{mass:padded-fillings}, this gives
$\Mass(h_j(\sigma))<\widehat b_j(s)$ and
$\dr h_j(\sigma)=z_\sigma$.
Choose zero on boundary simplices, and reuse each choice in every
incident simplex. The same cancellation for a $k$-simplex proves
$\dr T_f=0$, while its $k+1$ face terms give
\eqref{mass:initial-budget}.
\end{proof}

Above dimension $k$ we can no longer fill cycles, and individual fillings are replaced by families of small cycles supported in the simplices. These families play the role of the geometric chains $\iota_j(\sigma)$: the family attached to a simplex contracts the sum of the families attached to its faces.

\begin{lemma}\label{mass:simplex-families}
For the fixed $s>0$, the triangulation can be chosen fine enough that each $j$-simplex $\sigma$, $k+1\le j\le n$, has a mass-continuous family
\[
 c_j(\sigma):I^{j-k}\longrightarrow\ZC_k(M,\partial M;\Zt)
\]
supported in $\sigma$, with supremum of mass less than $s$. The families can be chosen with these face identities:
\begin{equation}\label{mass:simplex-base}
 c_{k+1}(\sigma)(0)=\dr\iota_{k+1}(\sigma),\quad
 c_{k+1}(\sigma)(1)=0;
\end{equation}
and, for $j\ge k+2$, $t\in I^{j-k-1}$, and $u\in I$,
\begin{equation}\label{mass:simplex-faces}
 c_j(\sigma)(t,0)=
 \sum_{\substack{\tau<\sigma\\\dim\tau=j-1}}c_{j-1}(\tau)(t),
\end{equation}
while $c_j(\sigma)$ is zero on all other boundary faces of its parameter cube. Boundary simplices have zero relative families. The map
\begin{equation}\label{mass:geometric-sphere}
 G(t)=\sum_{\sigma\in K_n}c_n(\sigma)(t),\quad t\in I^{n-k},
\end{equation}
vanishes on the parameter boundary and has Almgren class
\begin{equation}\label{mass:geometric-class}
 \Alm([G])=[M,\partial M].
\end{equation}
\end{lemma}
\begin{proof}
\noindent\textbf{Step 1: Construction of the families.}\enspace
Construct absolute-cycle families $\widehat c_j(\sigma)$ supported in $\sigma$. Set
$h=\max_{\sigma\in K_n}\operatorname{diam}_g(\sigma)$,
where diameter uses the Riemannian distance. Start with a finite
smooth triangulation and subdivide its fixed coordinate simplices
affinely, compatibly on faces. The original parametrizations and
their inverses restrict with uniform Lipschitz bounds; we do not
rescale each refined simplex to unit size. Thus $h\to0$ and
$\Mass(\partial\widehat\iota_{k+1}(\sigma))\le C h^k$.
Apply the supported contraction in Lemma~\ref{mass:small-cone}, in these coordinates, to obtain
\[
 \widehat c_{k+1}(\sigma)(0)=\partial\widehat\iota_{k+1}(\sigma),
 \quad \widehat c_{k+1}(\sigma)(1)=0.
\]

Inductively, for $j\ge k+2$, set
\[
 b_\sigma(t)=\sum_{\substack{\tau<\sigma\\\dim\tau=j-1}}
 \widehat c_{j-1}(\tau)(t),\quad t\in I^{j-k-1}.
\]
For $j=k+2$, the endpoints vanish because
\[
 b_\sigma(0)=\sum_{\substack{\tau<\sigma\\\dim\tau=k+1}}
 \partial\widehat\iota_{k+1}(\tau)
 =\partial^2\widehat\iota_{k+2}(\sigma)=0,
 \quad b_\sigma(1)=0.
\]
For $j\ge k+3$, the bottom face vanishes because
\[
 b_\sigma(t',0)=
 \sum_{\substack{\tau<\sigma\\\dim\tau=j-1}}
 \sum_{\substack{\rho<\tau\\\dim\rho=j-2}}
 \widehat c_{j-2}(\rho)(t')=0,
 \quad t'\in I^{j-k-2},
\]
each $\rho$ occurring twice. All other boundary values are zero termwise. Lemma~\ref{mass:small-cone} contracts $b_\sigma$ to zero within $\sigma$, fixing its parameter boundary. Choose the resulting $\widehat c_j(\sigma)$ once and reuse it in incident simplices.

\noindent\textbf{Step 2: Mass bounds and relative classes.}\enspace
With a common contraction constant $C$ for the finitely many parameter dimensions,
\[
 \begin{split}
 \sup_t\Mass(\widehat c_{k+1}(\sigma)(t))&\le C_{k+1}h^k,\\
 \sup_t\Mass(\widehat c_j(\sigma)(t))
 &\le C(j+1)\max_{\substack{\tau<\sigma\\\dim\tau=j-1}}
       \sup_t\Mass(\widehat c_{j-1}(\tau)(t))
 \le C_jh^k,
 \end{split}
\]
where $C_j=C(j+1)C_{j-1}$ is independent of $h$. Choose $h$ to satisfy every contraction threshold and $\max_j C_jh^k<s$.

The relative classes $c_j(\sigma)=[\widehat c_j(\sigma)]$ satisfy the stated identities and vanish on boundary simplices. For $n\ge k+2$, the bottom face of $G$ vanishes by
\[
 G(t,0)=\sum_{\sigma\in K_n}
 \sum_{\substack{\tau<\sigma\\\dim\tau=n-1}}c_{n-1}(\tau)(t)=0:
\]
interior faces occur twice and boundary-face families vanish. Other parameter faces are zero termwise. For $n=k+1$, use instead
\[
 G(0)=\dr\sum_{\sigma\in K_n}\iota_n(\sigma)=0,
 \quad G(1)=0.
\]

\noindent\textbf{Step 3: Computation of the Almgren class.}\enspace
We use Almgren's discrete chain construction \cite{Almgren1962}, in the relative mod-two formulation of \cite[Section 2.5]{LMN2018} and Theorem~\ref{inp:almgren}. Recall its principle: after a fine subdivision of the parameter cube, one fills the differences of the cycles at adjacent vertices by small chains, then fills the cycles formed by these chains around each two-dimensional cell, and so on; the sum of the chains attached to the top-dimensional cells represents the Almgren class. We carry this out for each family $\widehat c_j(\sigma)$, with fine parameter subdivisions matching on common faces, and show that the sum of its top-dimensional chains is the simplex $\widehat\iota_j(\sigma)$ itself. Inductively in $j$, assign to each $r$-cell $\alpha$ of the parameter subdivision for $\widehat c_j(\sigma)$ an absolute $(k+r)$-chain $T_\sigma(\alpha)$ supported in $\sigma$, satisfying
\begin{equation}\label{mass:parameter-chain-boundary}
 T_\sigma(v)=\widehat c_j(\sigma)(v),\quad
 \partial T_\sigma(\alpha)=
 \sum_{\substack{\beta<\alpha\\\dim\beta=r-1}}T_\sigma(\beta)
 \quad(r\ge1).
\end{equation}
For $j=k+1$, assign vertex values and fill adjacent differences by cones in the convex simplex coordinates. For $j\ge k+2$, first prescribe the parameter-boundary chains by
\begin{equation}\label{mass:parameter-chain-bottom}
 T_\sigma(\alpha\times\{0\})=
 \sum_{\substack{\tau<\sigma\\\dim\tau=j-1}}T_\tau(\alpha)
\end{equation}
on the bottom and zero on the other faces. Vertex values agree by Step~1, and \eqref{mass:parameter-chain-boundary} on bottom cells follows by summing the lower-dimensional identities. These prescriptions also agree at intersections: for $j=k+2$, the two bottom endpoints are $\partial^2\widehat\iota_{k+2}(\sigma)=0$ and zero; for $j\ge k+3$, the only non-termwise cancellation is
\[
 T_\sigma(\beta\times\{0\}\times\{0\})
 =\sum_{\substack{\tau<\sigma\\\dim\tau=j-1}}
   \sum_{\substack{\rho<\tau\\\dim\rho=j-2}}T_\rho(\beta)
 =0,
\]
where every codimension-two face occurs twice. Assign the remaining vertices by \eqref{mass:parameter-chain-boundary}, then fill remaining cells in increasing dimension by coning
\[
 Z_\alpha=\sum_{\substack{\beta<\alpha\\\dim\beta=r-1}}T_\sigma(\beta)
\]
inside $\sigma$. This is a cycle: for $r=1$ it is a difference of cycles, and for $r\ge2$ codimension-two parameter faces cancel. No prescribed chain is changed.
To check the smallness required by Almgren's construction, fix the
spatial triangulation first. On parameter edges, mass continuity
makes the boundary differences tend uniformly to zero, and coning
costs at most $C\operatorname{diam}(\sigma)$ times their mass.
Inductively, higher cell boundaries are sums of finitely many such
small chains, and the prescribed bottom chains obey the same estimate.
The relative chains $T_G(\alpha)=[\sum_{\sigma\in K_n}T_\sigma(\alpha)]$
agree with $G$ at vertices, satisfy \eqref{mass:parameter-chain-boundary}
in the relative quotient, and vanish on the parameter boundary by the
face cancellation of Step~2. For any $\eta>0$, refine the common
parameter grid until $\Mass(T_\sigma(\alpha))<\eta/\#K_n$ on every
positive-dimensional cell. Then $\Mass(T_G(\alpha))<\eta$ on each
such cell, so for sufficiently small $\eta$ the top-cell sum of $T_G$
computes $\Alm([G])$.

Define
\[
 Q_\sigma=
 \sum_{\substack{\alpha\text{ a cell of }I^{j-k}\\\dim\alpha=j-k}}
 T_\sigma(\alpha).
\]
Cancellation of interior parameter faces gives
\[
 \spt Q_\sigma\subset\sigma,\quad
 \partial Q_\sigma=
 \begin{cases}
  \partial\widehat\iota_{k+1}(\sigma),&j=k+1,\\
  \displaystyle\sum_{\substack{\tau<\sigma\\\dim\tau=j-1}}Q_\tau,&j\ge k+2.
 \end{cases}
\]
Inductively, therefore,
\[
 \partial\bigl(Q_\sigma-\widehat\iota_j(\sigma)\bigr)=0,
 \quad \spt\bigl(Q_\sigma-\widehat\iota_j(\sigma)\bigr)\subset\sigma.
\]
Pull the difference back to $\R^j$ and extend it by zero outside the simplex. The constancy theorem makes its mod-two multiplicity constant, hence zero by compact support. Thus $Q_\sigma=\widehat\iota_j(\sigma)$, and summing gives
\[
 \Alm([G])=
 \left[\sum_{\sigma\in K_n}Q_\sigma\right]
 =\left[\sum_{\sigma\in K_n}\widehat\iota_n(\sigma)\right]
 =[M,\partial M],
\]
as required.
\end{proof}

\begin{remark}\label{mass:concentration}
These families, and the finite sums and contractions used below, are mass-continuous on compact parameter spaces. By the argument following Definition~\ref{def:noconcentration}, they have no concentration of mass and satisfy the discretization hypotheses in Section~\ref{sec:inputs}.
\end{remark}

\begin{proof}[Proof of Theorem~\ref{thm:mass}]
Fix $N=N(n-k)$. We claim that for every $s>0$ there is a nonzero integral $k$-varifold $V_s$, stationary with free boundary and almost minimizing in small relative annuli, with
property $\mathcal P_N$ and $0<\|V_s\|(M)<L_n(s)$.
Suppose this assertion fails for some $s>0$. Fix a triangulation satisfying Lemmas~\ref{mass:chain-homotopy} and~\ref{mass:simplex-families}. We shall contract, simplex by simplex, all the families of cycles built from this triangulation, and reach a contradiction with $[M,\partial M]\ne0$.

Since $N(r)\le N$ for $r\le n-k$, we have $\mathcal P_{N(r)}\Rightarrow\mathcal P_N$. Thus alternative \textup{(i)} of Lemma~\ref{mass:sublevel} is excluded whenever $r\le n-k$ and $\sup\Mass(f)<L_n(s)$; every such family admits alternative \textup{(ii)}.
For this implication, note that the first $N(r)$ annuli of an admissible
collection of $N$ annuli again satisfy the separation inequalities,
with the same fixed radius $R_*$.

\noindent\textbf{Step 1: Initial contractions.}\enspace
For each $k$-simplex $f$, \eqref{mass:initial-budget} gives $\Mass(T_f)<L_k(s)-s<L_n(s)$. Lemma~\ref{mass:sublevel}\textup{(ii)}, with $r=0$ and $e=s$, gives a path
\[
 H_k(f):I\longrightarrow\ZC_k(M,\partial M;\Zt),\quad
 H_k(f)(0)=T_f,\quad H_k(f)(1)=0,
\]
with $\sup\Mass(H_k(f))<L_k(s)$. Set $H_k(f)=0$ when $f\subset\partial M$, and use the same chosen path in every simplex containing $f$.

\noindent\textbf{Step 2: Induction over the simplex dimension.}\enspace
Inductively, for $k+1\le j\le n$, define
\begin{equation}\label{mass:obstruction-sphere}
 f_j(\sigma)(t)=c_j(\sigma)(t)
 +\sum_{\substack{\tau<\sigma\\\dim\tau=j-1}}H_{j-1}(\tau)(t),
 \quad t\in I^{j-k},
\end{equation}
where each previously chosen $H_\ell$ has mass less than $L_\ell(s)$, has bottom value $f_\ell$ for $\ell>k$, and vanishes on the remaining parameter boundary. Choices are reused on common faces and zero on boundary simplices.

The family $f_j(\sigma)$ is based. For $j=k+1$, its initial value is
\[
 \dr\iota_{k+1}(\sigma)
  +\sum_{\substack{\tau<\sigma\\\dim\tau=k}}T_\tau=0,
\]
by cancellation of geometric faces and paired $h_{k-1}$ terms; its final value is zero.

For $j\ge k+2$, its bottom value is
\[
 \sum_{\substack{\tau<\sigma\\\dim\tau=j-1}}c_{j-1}(\tau)
 +\sum_{\substack{\tau<\sigma\\\dim\tau=j-1}}
   \left(c_{j-1}(\tau)
    +\sum_{\substack{\rho<\tau\\\dim\rho=j-2}}H_{j-2}(\rho)\right)=0.
\]
The $c_{j-1}$ terms cancel and each $H_{j-2}$ term occurs twice. Other parameter faces vanish termwise.

The $j+1$ faces give
\begin{equation}\label{mass:obstruction-budget}
 \sup\Mass(f_j(\sigma))
 <s+(j+1)L_{j-1}(s)=L_j(s)-s.
\end{equation}
For $j<n$, Lemma~\ref{mass:sublevel}\textup{(ii)} with $e=s$ therefore gives
a map $H_j(\sigma):I^{j-k+1}\to\ZC_k(M,\partial M;\Zt)$
whose bottom face is $f_j(\sigma)$ and whose other boundary faces are zero, with
\[
 \sup\Mass(H_j(\sigma))
 <\sup\Mass(f_j(\sigma))+s<L_j(s).
\]
Choose $H_j$ once per simplex and zero on boundary simplices, completing the induction. At $j=n$, retain the based families $f_n(\sigma)$ without further contraction.

\noindent\textbf{Step 3: The nonzero Almgren class and the contradiction.}\enspace
Equation~\eqref{mass:obstruction-sphere} gives
\[
 \sum_{\sigma\in K_n}f_n(\sigma)
 =\sum_{\sigma\in K_n}c_n(\sigma)=G.
\]
Indeed, every interior $(n-1)$-face appears twice, and every boundary-face family is zero. Additivity of the Almgren isomorphism and \eqref{mass:geometric-class} imply
\[
 \sum_{\sigma\in K_n}\Alm([f_n(\sigma)])
 =\Alm([G])=[M,\partial M]\ne0.
\]
Some $f_n(\sigma_*)$ therefore has nonzero Almgren class, necessarily $[M,\partial M]$ since $M$ is connected. Corollary~\ref{inp:sweepout} and \eqref{mass:obstruction-budget} give a nonzero stationary integral varifold $V_s$, almost minimizing in small relative annuli, with
\[
 0<\|V_s\|(M)\le\sup\Mass(f_n(\sigma_*))<L_n(s).
\]
It also satisfies $\mathcal P_N$, contradicting our assumption. Thus such a $V_s$ exists for every $s>0$.

\noindent\textbf{Step 4: Passage to the exact mass bound.}\enspace
Equation~\eqref{mass:budget-limit} gives \eqref{mass:limiting-bound}. To obtain the exact bound, set $s_\ell=1/\ell$ and choose $V_\ell=V_{s_\ell}$. They satisfy the same $\mathcal P_N$ and, by monotonicity,
\[
 0<\|V_\ell\|(M)<L_n(s_\ell)\le L_n(1).
\]
Lemma~\ref{inp:annular-compactness} gives a nonzero integral limit $V$, stationary with free boundary and almost minimizing in small relative annuli. Mass convergence and \eqref{mass:budget-limit} yield
\[
 0<\|V\|(M)=\lim_{\ell\to\infty}\|V_\ell\|(M)
 \le\lim_{\ell\to\infty}L_n(s_\ell)
 =\frac{(n+1)!}{k!}A.
\]
This proves \eqref{eq:mainmass}. In particular, the bound is $n(n+1)b_{n-2}$ for $k=n-1$.
\end{proof}

\begin{remark}\label{mass:NR}
The recursion and factorial coefficient follow the closed-manifold
construction of Nabutovsky--Rotman \cite[Theorem 1.4 and Section 4]{NR2006}.
The relative argument above uses mod-two cycles, begins with paths to
the boundary of length at most $D$, and retains compatible choices on
all parameter faces. The homology assumptions are used only to fill
individual cycles below dimension $k$; the final obstruction is
$[M,\partial M]\ne0$.
\end{remark}

Affine filling bounds make the recursion explicit.

\begin{corollary}
\label{hf:cor:recursion}
Under the homology hypotheses of Theorem~\ref{thm:mass},
suppose $\FH_j(L)\le\alpha_jL+\beta_j$ for
$1\le j\le k-1$, with $\alpha_j,\beta_j\ge0$.
Then its varifold can be chosen with
\begin{equation}\label{hf:eq:general-recursion}
 \|V\|(M)\le (n+1)!
 \left[
 D\prod_{j=1}^{k-1}\alpha_j+
 \sum_{i=1}^{k-1}\frac{\beta_i}{(i+1)!}
               \prod_{j=i+1}^{k-1}\alpha_j
 \right].
\end{equation}
An empty product is one.
\end{corollary}

\begin{proof}
In the recursion of Theorem~\ref{thm:mass}, bound
$b_1$ by $2\alpha_1D+\beta_1$ and each subsequent $b_j$
by $(j+1)\alpha_jb_{j-1}+\beta_j$. Dividing the $j$th bound by
$(j+1)!$ gives a recursion with multiplier $\alpha_j$ and additive
term $\beta_j/(j+1)!$, starting from
$b_1/2\le\alpha_1D+\beta_1/2$. Iteration and multiplication by
$(n+1)!$ give
\eqref{hf:eq:general-recursion}.
\end{proof}

\section{Closed filling estimates in dimension four}
\label{sec:closed-fill}
\label{sec:higher-fill}

The closed estimates use body--neck geometry in degrees one and two,
and BV theory and isoperimetry in degree three. Constants depend only
on the stated parameters; no numerical formula is asserted.
In degrees one and two the strategy is as follows. We cut the manifold into a bounded number of pieces of controlled geometry, called bodies and necks (Lemma~\ref{wz:lem:pieces}), and show that cycles which bound inside a piece admit fillings of linear cost in controlled ambient neighborhoods of the pieces (Lemmas~\ref{wz:lem:body} and~\ref{wz:lem:neck}). A given bounding cycle is then cut along the interfaces between the pieces, and the cuts are repaired at controlled cost (Propositions~\ref{wz:prop:closed} and~\ref{hf:prop:closed2}).

Let $\mathcal C(v_0,d_0)$ consist of smooth closed connected
four-manifolds $(X,h)$ satisfying
\begin{equation}\label{wz:eq:closed-class}
 |\operatorname{Ric}_h|\le3,\quad
 \Vol_h(X)\ge v_0>0,\quad \diam_h(X)\le d_0.
\end{equation}
For an absolute $j$-cycle $z$ in $X$ with coefficients in
$G\in\{\Zt,\mathbb Z\}$, put
\[
 \operatorname{Fill}_X(z)
 =\inf\{\Mass_h(Q):Q\in\II_{j+1}(X;G),\ \partial Q=z\},
\]
and define the filling function restricted to bounding cycles by
\begin{equation}\label{wz:eq:bounded-function}
 \operatorname{HF}^{\mathrm{bd}}_{j,G}(X;L)
 =\sup_{\substack{z\in\partial\II_{j+1}(X;G)\\\Mass_h(z)\le L}}
    \operatorname{Fill}_X(z).
\end{equation}
We write $\FH_j^{\mathrm{bd}}$ for the analogous function of relative
mod-two cycles of $(M,\partial M)$ that bound. These functions differ
from the all-cycle functions when homology does not vanish; in
particular, $\FH_j^{\mathrm{bd}}=\FH_j$ if $H_j(M,\partial M;\Zt)=0$.

Wu--Zhu \cite{WZ2020} prove an integer all-cycle estimate in degree
one under $H_1(X;\mathbb Z)=0$. Fu--Zhu
\cite[Theorems 1.1--1.2]{FZ2026} treat the Einstein subclass in
that degree, with quantitative analytic and combinatorial inputs.
Neither supplies the degree-two or relative statements proved below.

Our geometric input, also used by Wu--Zhu, is the body--neck
construction of Cheeger--Naber \cite[Section 8]{CN2015}.

\begin{lemma}\label{wz:lem:pieces}
A manifold in $\mathcal C(v_0,d_0)$ can be cut into a bounded
number of compact body and neck pieces with the following properties.
Here all constants and the bound on the number of pieces depend
only on $v_0,d_0$.
\begin{enumerate}[label=\textup{(\roman*)},leftmargin=2.2em]
\item A body $K$, including its cutting cross-sections, has a
scale $s\le C d_0$, ambient diameter at most $Cs$, and harmonic radius
at least $cs$ throughout the ambient $cs$-neighborhood of $K$,
after decreasing $c$ if necessary.
\item Each neck is uniformly bilipschitz to a cone annulus
\[
 [r,R]\times Y,\quad d\rho^2+\rho^2g_Y,\quad
 Y=S^3/\Gamma,\quad |\Gamma|\le N,
\]
where $0<r<R\le C d_0$. No bound on $R/r$ is required.
\item The interfaces are disjoint spherical-space-form
cross-sections. Each interface can be moved through a buffered
collar of width comparable to its radius, retaining \textup{(i)}
and \textup{(ii)}. The collars have uniformly bounded multiplicity.
\end{enumerate}
\end{lemma}

\begin{proof}
After rescaling by a factor depending only on $v_0,d_0$,
\cite[proof of Theorem 8.63]{CN2015} applies. That theorem, Theorem 8.6
and Lemma 8.39 give the cone, group-order and harmonic-radius
estimates. The decomposition constructed there consists of overlapping
open sets, and we have to trim these sets into pieces with disjoint
interfaces. We describe the cuts at one scale, then the root of the
decomposition, and finally the recursion.

At each parent scale, \cite[proof of Lemma 8.60]{CN2015} selects
good radii $r_i$ bounded below by a fixed multiple of that scale,
with disjoint balls $B_{r_i/4}(x_i^0)$. Let $\Gamma_i$ be the group
of the cone at this good scale. If $\Gamma_i$ is trivial, the
harmonic radius is at least $c r_i$ throughout $B_{3r_i/4}(x_i^0)$.
Indeed, rescale $r_i$ to one. As the good-scale tolerance and
the rescaled Ricci bound tend to zero, the compactness argument
in \cite[Theorem~8.6, (8.8)--(8.10)]{CN2015} gives a Euclidean
limit ball $B_2(0)\subset\R^4$ in this case. If the claimed bound
failed, a subsequence of offending points in $B_{3/4}(x_i^0)$ would
therefore have balls of radius $1/2$ converging to a Euclidean ball
of that radius. The rescaled manifolds are uniformly noncollapsed,
so \cite[Theorem~2.11]{CN2015} gives harmonic radius at least $1/4$
at these points for all sufficiently large indices. Fix the tolerance
accordingly. The half-ball is regular at the body scale, and we keep
its disjoint quarter-ball in the body and stop this branch.

For each nontrivial $\Gamma_i$, apply the center adjustment of
Lemma~8.60 to obtain $d(x_i,x_i^0)\le\eps r_i$.
Fix $\eps\ll q\ll1/4$ and cut
at cone-coordinate radii in $[qr_i,2qr_i]$. After decreasing $q$,
Lemma 8.39 places these collars in disjoint balls $B_{Cqr_i}(x_i)$:
the original disjoint balls give separation at least $(r_i+r_j)/4$,
moving the centers loses at most $\eps(r_i+r_j)$, and we choose
$Cq<1/4-\eps$. The proof of Lemma 8.60 also gives harmonic-radius
control throughout $A_{\eps r_i,2r_i}(x_i)$ and places the inner end of
the neck inside $B_{\eps r_i}(x_i^0)$. Thus the proposed cuts lie in the
cone chart, with a collar of width comparable to $qr_i$ on both sides.
The interfaces are smooth coordinate links, not distance spheres.
We absorb the outer neck portions into the body; their harmonic radius
is at least $c(q)r_i$, which preserves the body-scale estimate.

At the root, volume comparison gives local noncollapse, and
\cite[Lemma 8.15 and Corollary 8.21]{CN2015} give good dyadic radii
in a common finite range $[r_*,r^*]$. Select disjoint quarter-balls
maximally in decreasing order of radius. Each omitted center lies
within a selected half-ball, because its quarter-ball meets a selected
one of at least its radius. These half-balls therefore cover $X$,
and their number is bounded by packing. Outside the $q$-scale holes
each point belongs either to a terminal half-ball or to the regular
annulus of a selected ball with nontrivial group.
Since $r_i\ge r_*$, the enlarged root body and a fixed-width
neighborhood of its cuts have uniform harmonic radius.

For the inner cut of this neck, write $x=x_i$ and let $s=r_\beta$
be its stopping radius in \cite[Lemma 8.60]{CN2015}.
The inclusion $B_{2s}(x)\subset B_{\eps r_i}(x_i^0)$ in the proof of
that lemma gives $s\le\eps r_i$. Thus the inner cut described below
and the outer cuts at radii $qr_i$ to $2qr_i$ have disjoint collars.
Let $\rho$ denote the cone radius. On the metric annulus
$A_{3s/4,\,3s/2}(x)$ we may take
\[
  |\rho(y)-d(x,y)|<s/32.
\]
Indeed, after rescaling by $s$, the dyadic charts are obtained from
pointed convergence to the cone with $x$ converging to its vertex
\cite[proof of Theorem 8.6 and (8.45)--(8.46)]{CN2015}.
On the indicated compact subannulus their radial coordinates therefore
converge uniformly to distance from $x$ as the tolerance tends to zero.
The rotations used to align the charts preserve radius, and the gluing
changes the inverse charts there by at most $C\eps s$
\cite[(8.48)--(8.52)]{CN2015}. Fix the tolerance so that the displayed
estimate holds.

Cut at any level $9s/8\le t\le5s/4$. The enclosed region is
\[
  B_{7s/8}(x)\ \cup\
  \{y\in A_{3s/4,\,3s/2}(x):\rho(y)\le t\}.
\]
On the overlap with the inner ball one has $\rho<t$, and near the
outer edge of the annulus one has $\rho>t$. The level in this annulus
is nonempty by continuity along a minimizing geodesic from $x$ to its
outer edge. It is both open and closed in the connected coordinate
link $\rho=t$, so the boundary is precisely that link.
The region contains $B_s(x)$ and
lies in $B_{4s/3}(x)$. The entire range of cuts has a collar on both
sides in the same chart. Outside $B_{7s/8}(x)$ the enclosed region
already has harmonic radius at least $cs$ by
\cite[Theorem 8.6]{CN2015}.

It remains to cover $\overline B_{7s/8}(x)$. Choose the initial
scale for the good-radius selection to be the smaller of $s/100$
and the fixed multiple of $s$ in \cite[(8.73)]{CN2015}.
The dyadic selection of Lemma 8.15 then gives a candidate radius
between $c s$ and this initial scale at every point of the core.
As at the root, take a maximal disjoint collection of quarter-balls
in decreasing order of radius. The half-balls cover the core, and
even after recentering the selected balls of twice the candidate
radius lie in $B_s(x)$. In particular all new cut collars are
separated from the incoming cuts. Together the regular parts of
these half-balls and the incoming collar give the child body, of
ambient diameter at most $3s$ and harmonic radius at least $cs$. The same
bound holds on a smaller $cs$-neighborhood by the buffers at the cuts
and the local defining property of the harmonic radius.

All child centers lie in $B_{2s}(x)$, where the volume-ratio bound
in \cite[(8.73)]{CN2015} holds. Bishop--Gromov monotonicity preserves
the bound at the smaller candidate radii, so \cite[(8.75)]{CN2015}
gives a strict decrease in group order at every nontrivial branch.
This bounds the depth, and packing bounds the number of children
at each stage. The number of pieces and the collar multiplicity
are therefore uniformly bounded. Lemma 8.39 controls each trimmed
neck independently of $R/r$.
\end{proof}

\begin{lemma}\label{wz:lem:body}
For a body $K$ of scale $s$ in Lemma~\ref{wz:lem:pieces},
let $j=1$ or $2$ and $G=\Zt$ or $\mathbb Z$.
Every $j$-cycle $z$ with coefficients in $G$ that bounds a
$(j+1)$-current supported in $K$ has a filling in a controlled ambient neighborhood of
$K$ satisfying
\begin{equation}\label{wz:eq:body}
 \operatorname{Fill}_X(z)\le Cs\,\Mass_h(z).
\end{equation}
The constant is uniform over $\mathcal C(v_0,d_0)$.
\end{lemma}

\begin{proof}
By scaling we may assume $s=1$; the rescaled Ricci bounds remain
uniform since $s\le Cd_0$. The idea is to transport the problem, by a
bilipschitz map with uniform constants, to one of finitely many fixed
reference domains. There, a linear filling inequality follows from the
deformation theorem and linear algebra.

We first produce the reference domains.
Harmonic-radius control, bounded diameter and volume comparison give
uniform covering numbers, and hence marked Gromov--Hausdorff
precompactness, for closed neighborhoods marked by $K$. Here the mark
records the compact subset $K$, not only its ambient neighborhood.
Covering numbers are uniform because harmonic balls have a definite
volume and the enclosing ball has uniformly bounded volume; the same
finite covers give precompactness of the marked subsets in the
Hausdorff metric. Choose a finite fine marked net of reference pairs.
The construction in \cite[Theorem 8.1]{CN2015} gives uniformly
bilipschitz diffeomorphisms $\Phi$ between buffered neighborhoods of $K$
and of its reference $K_\alpha$. For $\rho$ below the harmonic scale,
we arrange
\begin{equation}\label{wz:eq:marked-buffers}
 \Phi(K)\subset T_{\rho/8}(K_\alpha),\quad
 T_{\rho/2}(K_\alpha)\subset\operatorname{image}\Phi,
\end{equation}
where $T_t$ denotes the ambient $t$-neighborhood.

The inclusions \eqref{wz:eq:marked-buffers} are obtained as follows.
To track the core $K$, start with a marked Gromov--Hausdorff
approximation $f$ of error $\delta$. Harmonic-coordinate smoothing at
radius $\sigma$ has displacement at most $C(\sigma+\delta)$: points
within $\sigma$ have images within $\sigma+2\delta$, and averaging in
bilipschitz charts preserves this bound; bounded-overlap patching
changes only $C$. Thus $\sup d(\Phi,f)\le C(\sigma+\delta)$ on a smaller
neighborhood. Choose $\sigma$ small relative to $\rho$, and then
$\delta$ small enough for this tracking and for the diffeomorphism and
metric comparison of Theorem 8.1. The marked approximation gives the
first inclusion in \eqref{wz:eq:marked-buffers}, and the buffered
neighborhood inclusions of that theorem give the second.

Next we place the cycle in a fixed domain.
For each reference choose a fixed compact smooth domain
$L_\alpha$ with
\[
 T_{\rho/8}(K_\alpha)\subset\Int L_\alpha
 \subset L_\alpha\subset T_{\rho/4}(K_\alpha).
\]
Then $\Phi^{-1}$ is defined on $L_\alpha$, and $L_\alpha$ contains the
pushforward of any certificate $\partial Q=z$ supported in $K$.
In particular, $\Phi_\#z$ bounds in $L_\alpha$, even though the
efficient filling constructed below need not stay in $\Phi(K)$.

We now prove a linear filling inequality in $L_\alpha$.
Fix a finite Lipschitz triangulation of $L_\alpha$ in which the
boundary is a subcomplex. The deformation theorem \cite{Federer1969},
applied in this triangulated domain, replaces a bounding
$j$-cycle $z'$ by a simplicial cycle $c$ and an error filling
$H$ such that
\[
 z'=c+\partial H,\quad
 \Mass(H)+\|c\|_{\ell^1}\le C_\alpha\Mass(z').
\]
Since $z'$ bounds in $L_\alpha$, so does $c$ as a current. By the
identification of current homology with simplicial homology, or by the
same deformation applied to a qualitative filling of $c$, the cycle $c$
lies in the simplicial image $\operatorname{im}\partial_{j+1}$; no mass
bound for that qualitative filling is needed. It remains to find a
simplicial chain $d$ with $c=\partial d$ and
$\|d\|_{\ell^1}\le C_\alpha\|c\|_{\ell^1}$. Over $\mathbb Z$, the group
$\operatorname{im}\partial_{j+1}$ is free abelian; choose a basis and
fixed lifts of the basis elements. The inverse of the basis map extends
to the real span of this group, and by equivalence of norms on a fixed
finite-dimensional space the basis coordinates of $c$ are bounded by
its $\ell^1$ norm. The chosen lifts then give $d$ with the stated bound.
Over $\Zt$, the group $\operatorname{im}\partial_{j+1}$ is finite:
choose one lift for each nonzero boundary, and use that every nonzero
simplicial chain has $\ell^1$ norm at least one. In both cases $H+d$
fills $z'$ at cost $C_\alpha\Mass(z')$.

To conclude, apply this inequality to $z'=\Phi_\#z$, take the maximum of the constants over
the finitely many references, pull the filling back by $\Phi$, and
rescale; the rescaling produces the factor $s$ in \eqref{wz:eq:body}.
\end{proof}

The efficient filling may leave the cut body but remains in $X$.
Finite diffeomorphism types alone would not suffice: the marked
comparison also places an existing filling in the fixed domain.

\Needspace{7\baselineskip}
\begin{lemma}\label{wz:lem:neck}
There is a uniform constant $C$ with the following properties.
On a radius-$r$ link $L$ of Lemma~\ref{wz:lem:pieces}:
\begin{enumerate}[label=\textup{(\roman*)},leftmargin=2.2em]
\item a zero-cycle $\gamma$ that bounds on $L$ has a
one-dimensional filling of mass at most $Cr\,\Mass(\gamma)$;
\item each class in $H_1(L;\Zt)$ has a representative of
mass at most $Cr$;
\item a one-cycle that bounds on $L$ has a two-dimensional
filling of mass at most $Cr$ times its length.
\end{enumerate}
A one-cycle $z$ that bounds in a neck with outer radius $R$
has a filling there with mass at most $CR\,\Mass(z)$.
\end{lemma}

\begin{proof}
Since $|\Gamma|\le N$, there are finitely many groups $\Gamma$ and
finitely many four-dimensional orthogonal representation types, hence
finitely many isometry types of round space forms $S^3/\Gamma$. It
therefore suffices to prove \textup{(i)--(iii)} for one fixed space
form of unit size; scaling by $r$ and the uniform bilipschitz
comparison of Lemma~\ref{wz:lem:pieces} then give the statements for
the links. At unit scale, we fill a bounding zero-cycle by pairing its
points and joining each pair by a minimizing geodesic, we choose fixed
representatives of the finitely many mod-two homology classes, and we
fill bounding one-cycles using a fixed triangulation, the deformation
theorem, and the boundary-matrix argument in the proof of
Lemma~\ref{wz:lem:body}.

For the last assertion, consider first the exact cone
$[r,R]\times Y$, and let $p(\rho,y)=(r,y)$ be the projection onto the
inner link. The radial homotopy
$(\rho,y)\mapsto((1-u)\rho+ur,y)$, $0\le u\le1$, from the identity to
$p$ is spatially $1$-Lipschitz and has speed at most $R$, so the
homotopy formula fills $z-p_\#z$ at cost $R\Mass(z)$. Since $z$ bounds
in the neck, $p_\#z$ bounds on the inner link, and \textup{(iii)} fills
it at cost $Cr\Mass(p_\#z)\le Cr\Mass(z)$. The sum of the two fillings
has mass at most $CR\Mass(z)$.
Bilipschitz comparison gives the neck estimate independently of $R/r$.
\end{proof}

We now prove a mod-two bounding-cycle counterpart of
\cite[Theorem 1.6]{WZ2020}, keeping track of homology on each piece.

\begin{proposition}
\label{wz:prop:closed}
There are $a(v_0,d_0),b(v_0,d_0)>0$ such that every
$(X,h)\in\mathcal C(v_0,d_0)$ satisfies
\begin{equation}\label{wz:eq:closed-bound}
 \operatorname{HF}^{\mathrm{bd}}_{1,\Zt}(X;L)
 \le a(v_0,d_0)L+b(v_0,d_0).
\end{equation}
No vanishing of absolute homology is assumed.
\end{proposition}

\begin{proof}
Let $z$ be a bounding one-cycle, and fix $Q$ with $\partial Q=z$. The
chain $Q$ is used only to certify that certain cycles bound; its mass
plays no role. We cut $z$ along the interfaces of
Lemma~\ref{wz:lem:pieces}, close up the resulting pieces by short caps,
correct the caps so that each piece bounds inside its own piece of the
manifold, and fill the pieces by Lemmas~\ref{wz:lem:body}
and~\ref{wz:lem:neck}.

We begin with the cuts.
In each interface collar, of scale $r_e$, the coarea formula gives a
level that is regular for both $z$ and $Q$, carries no mass of either,
and on which the trace $\gamma_e$ of $z$ satisfies
\begin{equation}\label{wz:eq:trace-bound}
 \Mass(\gamma_e)\le
 \frac{C}{r_e}\Mass(z\res\text{collar}_e),
\end{equation}
while the slice $Q_e$ of $Q$ satisfies $\partial Q_e=\gamma_e$.
To find such a level, integrate the slice mass of
$z$ over a collar interval of length comparable to $r_e$.
Markov's inequality leaves a set of levels of positive measure on which
\eqref{wz:eq:trace-bound} holds after doubling $C$; intersect it
with the full-measure sets of levels at which the slices of $z$ and $Q$
are integral and which are charged by neither mass measure. There are
only finitely many interfaces, so all cuts can be made with these
properties.

Next we cap the cut cycle.
Since $\gamma_e=\partial Q_e$ bounds on the link,
Lemma~\ref{wz:lem:neck}(i) and \eqref{wz:eq:trace-bound} provide a cap
$P_e^0$ on this link with $\partial P_e^0=\gamma_e$ and
$\Mass(P_e^0)\le C\Mass(z\res\text{collar}_e)$. The cycle $P_e^0-Q_e$
may represent a nontrivial class in the first homology of the link.
We correct it by a representative $R_e$ of this class from
Lemma~\ref{wz:lem:neck}(ii):
\[
 P_e=P_e^0-R_e,\quad
 \Mass(R_e)\le Cr_e,\quad P_e-Q_e=\partial T_e.
\]
Here $T_e$ is supported on the link; its mass is irrelevant.

It remains to fill the pieces.
For each piece $U$, cap $z\res U$ by the $P_e$ to obtain
$z_U=z\res U+\sum_{e\subset\partial U}P_e$. Modulo two,
$\partial(Q\res U)=z\res U+\sum_{e\subset\partial U}Q_e$, and
consequently $Q\res U+\sum_{e\subset\partial U}T_e$ has boundary $z_U$.
This verifies the local bounding hypothesis required by the body and
neck lemmas, without using the mass of $Q$. Caps cancel between
adjacent pieces, so $\sum_Uz_U=z$, and bounded collar multiplicity gives
\begin{equation}\label{wz:eq:piece-sum}
 \sum_U\Mass(z_U)\le C\Mass(z)+C\sum_e r_e.
\end{equation}
Lemmas~\ref{wz:lem:body} and~\ref{wz:lem:neck} fill each $z_U$ with
coefficient at most $Cd_0$. Since the interfaces are bounded in number
and $r_e\le Cd_0$, the sum of these fillings has boundary $z$ and mass
at most
\[
 C(v_0,d_0)\bigl(d_0\Mass(z)+d_0^2\bigr).
\]
This proves \eqref{wz:eq:closed-bound}; overlap of fillings is harmless.
\end{proof}

\begin{remark}\label{wz:rem:charges}
The additive term pays for $R_e$: short caps alone need not make
the piece cycles bound. If $H_1(X;\Zt)=0$,
\eqref{wz:eq:closed-bound} applies to all one-cycles; otherwise
it applies only to those that bound.
\end{remark}

The same decomposition controls bounding two-cycles; the link
homology distinguishes integer from mod-two coefficients.

\begin{proposition}\label{hf:prop:closed2}
There are positive constants $a_2,\beta_2,A_2$, depending only on
$v_0,d_0$, such that every
$X^4\in\mathcal C(v_0,d_0)$ satisfies
\begin{equation}\label{hf:eq:closed2}
 \operatorname{HF}^{\mathrm{bd}}_{2,\Zt}(X;L)
 \le a_2L+\beta_2.
\end{equation}
For oriented $X$ there is also an integer-coefficient estimate
\begin{equation}\label{hf:eq:closed2int}
 \operatorname{HF}^{\mathrm{bd}}_{2,\mathbb Z}(X;L)
 \le A_2L.
\end{equation}
The proof gives bounds of the forms
$C(v_0,d_0)(d_0L+d_0^3)$ and $C(v_0,d_0)d_0L$,
respectively. These are all-cycle estimates if the corresponding
$H_2(X;G)$ vanishes.
\end{proposition}

\begin{proof}
Use the pieces of Lemma~\ref{wz:lem:pieces}. On their radius-$r$
links, finite triangulations give linear filling constants $Cr$
in degrees one and two for either coefficient group. The finitely
many mod-two $H_2$ classes have representatives of mass at most $Cr^2$.

For integer coefficients,
\begin{equation}\label{hf:eq:link-h2}
 H_2(S^3/\Gamma;\mathbb Z)=0.
\end{equation}
Indeed spherical three-space forms are orientable with finite
integer $H_1$, so duality gives $H_2=H^1=0$. In contrast,
$H_2(\mathbb{RP}^3;\Zt)=\Zt$.

As in Lemma~\ref{wz:lem:neck}, radial contraction costs $R\Mass(z)$
for two-cycles, and the inner-link projection of a bounding cycle
fills at cost $Cr\Mass(z)$. Thus the neck constant is $CR$,
independently of $R/r$; Lemma~\ref{wz:lem:body} controls bodies.

Fix $z=\partial Q$ with $Q$ three-dimensional. Choose common
regular cuts for $z,Q$, carrying no mass of either, as in
Proposition~\ref{wz:prop:closed}. With consistent slice signs,
\[
 \partial Q_e=\gamma_e,\quad
 \Mass(\gamma_e)\le\frac{C}{r_e}
       \Mass(z\res\mathrm{collar}_e).
\]
Choose link caps $P_e^0$ with $\partial P_e^0=\gamma_e$ and
$\Mass(P_e^0)\le C\Mass(z\res\mathrm{collar}_e)$.
Over $\mathbb Z$, set $R_e=0$ using \eqref{hf:eq:link-h2};
over $\Zt$, choose $R_e$ representing $[P_e^0-Q_e]$.
In both cases set $P_e=P_e^0-R_e$, so that
$\Mass(R_e)\le Cr_e^2$ and $P_e-Q_e=\partial T_e$
for a link three-current $T_e$. Neither $Q$ nor $T_e$ is estimated.

Choose the orientation of each link once and write
$\partial(Q\res U)=z\res U+\sum_{e\subset\partial U}\epsilon_{Ue}Q_e$,
where the signs $\epsilon_{Ue}=\pm1$ are opposite for the two pieces
adjacent to $e$; modulo two all signs are omitted. Put
\[
 z_U=z\res U+\sum_{e\subset\partial U}\epsilon_{Ue}P_e,
 \quad
 Q_U=Q\res U+\sum_{e\subset\partial U}\epsilon_{Ue}T_e.
\]
Then $\partial Q_U=z_U$, which certifies that $z_U$ bounds in $U$, and
the caps cancel in $\sum_Uz_U=z$. Hence
\[
 \sum_U\Mass(z_U)\le
 \begin{cases}
 C\Mass(z),&G=\mathbb Z,\\
 C\Mass(z)+C\sum_e r_e^2,&G=\Zt.
 \end{cases}
\]
The body and neck coefficients are at most $Cd_0$; the interfaces
are bounded in number and $r_e\le Cd_0$. Summing the fillings of the
cycles $z_U$ therefore proves both assertions.
\end{proof}

\begin{remark}\label{hf:rem:coeff}
The mod-two additive term vanishes in this proof if every link
has $H_2(S^3/\Gamma;\Zt)=0$; odd group order suffices. We do not
claim necessity of that term in general. Nor does
\eqref{hf:eq:closed2int} imply \eqref{hf:eq:closed2}: mod-two surface
cycles need not admit mass-preserving integral lifts.

Wu--Zhu's hypothesis $H_1(X;\mathbb Z)=0$ alone cannot give an
all-cycle $HF_2$ bound: a factor sphere in $S^2\times S^2$ does not bound.
\end{remark}

In top degree, BV theory replaces the body--neck argument.
Let $M^n$ be smooth, compact and connected, possibly closed.
For a set of finite perimeter $E$ put
$P_M(E)=P(E;\Int M)$, taking $\Int M=M$ in the closed case.
Define the relative Cheeger constant by
\[
 h_N(M)=
 \inf_{0<\Vol(E)<\Vol(M)}
 \frac{P_M(E)}
      {\min\{\Vol(E),\Vol(M\setminus E)\}}.
\]
For closed $M=X$ we write $h(X)$.

\begin{proposition}\label{hf:prop:top}
Every bounding relative mod-two $(n-1)$-cycle $z$ satisfies
\begin{equation}\label{hf:eq:top}
 \inf_{\dr Q=z}\Mass(Q)
 \le\min\left\{\frac{\Vol(M)}2,\ h_N(M)^{-1}\Mass(z)\right\}.
\end{equation}
The corresponding closed statement holds with $h(X)$.
For oriented manifolds the integer-coefficient version
has the linear bound $h_N(M)^{-1}\Mass(z)$.
\end{proposition}

\begin{proof}
A mod-two filling is $\llbracket E\rrbracket$ for a finite-perimeter
set. By relative constancy, its only alternative is
$\llbracket M\setminus E\rrbracket$. Since $P_M(E)=\Mass(z)$,
the definition of $h_N(M)$ proves \eqref{hf:eq:top}.

For integer coefficients on oriented $M$, write a filling as
$u\llbracket M\rrbracket$, where $u$ is integer-valued and BV,
and choose an integer median $m$: both $\{u>m\}$ and $\{u<m\}$
have volume at most $\Vol(M)/2$. Every positive superlevel of
$u-m$ and every positive superlevel of $m-u$ therefore has volume
at most half. Apply the defining Cheeger inequality to these
level sets and integrate. The BV coarea formula
\cite{EG1992,Simon2014}, applied to the positive and negative
levels of $u-m$, gives
\[
 \Mass((u-m)\llbracket M\rrbracket)
 =\int_M|u-m|
 \le h_N(M)^{-1}|Du|(\Int M)
 =h_N(M)^{-1}\Mass(z).
\]
Subtracting $m\llbracket M\rrbracket$ preserves the relative boundary.
\end{proof}

The half-volume bound is specific to mod-two coefficients.
For example, over $\mathbb Z$ the cycle
$N\partial\llbracket E\rrbracket$ requires filling mass
$N\min\{\Vol(E),\Vol(M\setminus E)\}$.

\begin{lemma}\label{hf:lem:cheeger}
If a closed connected $X^n$ satisfies
$\operatorname{Ric}_X\ge-(n-1)\kappa$,
$\kappa\ge0$, and $\diam(X)\le d$, then
\begin{equation}\label{hf:eq:cheeger}
 h(X)^{-1}\le\int_0^{d/2}\cosh^{n-1}(\sqrt\kappa\,t)\,dt.
\end{equation}
\end{lemma}

\begin{proof}
This is Gallot's isoperimetric theorem \cite[Theorem 6.14]{Gallot1988},
with curvature parameter $K=-\kappa$. Gallot's normalization in
Section~6.C is perimeter divided by the smaller enclosed volume; smooth
BV approximation gives the finite-perimeter formulation used here.
\end{proof}

For an alternative route to Ricci--diameter control, see the
segment inequality in \cite[Theorem 2.11]{CC1996}. A manifold
$X^4\in\mathcal C(v_0,d_0)$ satisfies $\operatorname{Ric}_h\ge-3$, so
Proposition~\ref{hf:prop:top} and Lemma~\ref{hf:lem:cheeger} with
$\kappa=1$ yield
\begin{equation}\label{hf:eq:closed3}
 \operatorname{HF}^{\mathrm{bd}}_{3,\Zt}(X;L)
 \le\min\{\Vol(X)/2,c_3(d_0)L\},
\end{equation}
where
\begin{equation}\label{hf:eq:c3}
 c_3(d_0)=\int_0^{d_0/2}\cosh^3t\,dt
 =\sinh(d_0/2)+\tfrac13\sinh^3(d_0/2).
\end{equation}
Only lower Ricci and diameter bounds enter here. These do not control
$h_N(M)^{-1}$ for domains: smooth Euclidean dumbbells with fixed-size
lobes and shrinking tubes have bounded diameter, but cuts of vanishing
area separate fixed positive volumes. Such domains can be topological
balls. The relative half-volume bound still holds.

The body--neck argument relies on four-dimensional geometry and
spherical three-links; it gives no arbitrary-dimensional intermediate
degree estimate. The top-degree BV argument and the following transfer
are dimension independent.

\section{From closed to relative filling}
\label{sec:relative-wz}

We transfer the closed estimates by doubling and cutting in every
degree, or by completion through a complementary region in degree one.

The doubling argument is independent of Ricci curvature.
Write $X=DM=M_+\cup_S M_-$ for a smooth double, with
$S=\partial M$ and identifications $j_\pm$.

\begin{lemma}\label{wz:lem:double-chain}
For each $j\ge0$ the formula
\[
 \mathcal D[T]=(j_+)_\#T+(j_-)_\#T
\]
defines a chain map from relative integral mod-two chains
of $(M,\partial M)$ to absolute mod-two chains of $X$.
If $\Lambda^{-2}g\le j_\pm^*h\le\Lambda^2g$, then
\begin{equation}\label{wz:eq:double-mass}
 \Mass_h(\mathcal Dz)\le2\Lambda^j\Mass_g(z).
\end{equation}
In particular, a relative boundary doubles to an absolute
boundary, whether or not $H_j(X;\Zt)$ vanishes.
\end{lemma}

\begin{proof}
Boundary-supported chains cancel between the halves modulo two,
so the map is well-defined and commutes with boundary. For each
integral representative $T$,
$\Mass_h(\mathcal D[T])\le2\Lambda^j\Mass_g(T)$.
Taking the infimum proves \eqref{wz:eq:double-mass}, and
$\dr Q=z$ implies $\partial\mathcal DQ=\mathcal Dz$.
\end{proof}

\begin{lemma}\label{wz:lem:cut}
Let $N$ be a smooth compact domain in a closed Riemannian
manifold $(X,h)$ of dimension $n$, let $1\le j\le n-1$,
and let $B$ be an integral mod-two $(j+1)$-current
whose boundary agrees in $\Int N$ with an integral
representative $T$ of a relative $j$-cycle $z$ in $N$.
For a sequence of regular levels $t\downarrow0$ there are
relative fillings $Q_t$ of $z$ such that
\begin{equation}\label{wz:eq:cut}
 \Mass_h(Q_t)\le
 (1+Ct)^{j+1}\Mass_h(B\res N_t)+
 t(1+Ct)^j\Mass_h(T).
\end{equation}
Here $N_t$ is obtained by removing a collar of width $t$.
For two halves separated by their common boundary, the levels
can be chosen simultaneously, and the masses of the two
restricted portions of $B$ sum to at most $\Mass_h(B)$.
\end{lemma}

\begin{proof}
The difficulty is that $B$ may meet $\partial N$ irregularly, so that
$B\res N$ need not be an integral current. We therefore restrict $B$ to
the slightly smaller domain $N_t$ and compress the remaining collar onto
$\partial N$.

We first construct the compression.
In inward normal coordinates $(y,s)$ near $\partial N$, fix a smooth
cutoff $0\le\chi\le1$, equal to one near $s=0$ and supported in a
collar. For small $t>0$, the map
\[
 f_t(y,s)=(y,\max\{0,s-t\chi(s)\})
\]
fixes $\partial N$ pointwise, is the identity away from the collar,
and has Lipschitz constant at most $1+Ct$. Moreover, once $t$ is so
small that $\chi=1$ on $[0,t]$, it maps the collar
$N\setminus N_t=\{s\le t\}$ into $\partial N$.
The homotopy $u\mapsto f_{ut}$, $0\le u\le1$, from the identity to
$f_t$ has speed at most $t$. Since it is stationary on $\partial N$,
which contains the support of $\partial T$, the homotopy formula
supplies a current $E_t$ with
\[
 \dr E_t=z-[(f_t)_\#T],\quad
 \Mass_h(E_t)\le t(1+Ct)^j\Mass_h(T).
\]

Next we restrict the filling.
For almost every small $t$, the restriction $B_t=B\res N_t$ is an
integral current, and the restriction formula gives
$\partial B_t=T\res N_t+\Sigma_t$, where $\Sigma_t$ is the slice of $B$
at the level $s=t$; no restriction of $B$ to the original seam is
asserted to be integral. The map $f_t$ sends the slice $\Sigma_t$, and
also the omitted portion $T\res(N\setminus N_t)$ of $T$, into
$\partial N$, so both vanish in the relative quotient. Hence
$\dr[(f_t)_\#B_t]=[(f_t)_\#T]$.

Finally, set $Q_t=[(f_t)_\#B_t]+E_t$. Then $\dr Q_t=z$, and the Lipschitz bound
for $f_t$ together with the estimate for $E_t$ gives
\eqref{wz:eq:cut}. For two halves, use regular levels common to both.
The two regions $N_t$ are disjoint, so the masses of the two restricted
portions of $B$ sum to at most $\Mass_h(B)$; the mass of $B$ near the
seam is discarded.
\end{proof}

\begin{proposition}
\label{wz:prop:transfer}
Let $M^n$ have a smooth double $X$ with
$\Lambda^{-2}g\le j_\pm^*h\le\Lambda^2g$, and let
$1\le j\le n-1$. Suppose the closed mod-two bounding-cycle
function obeys
$\operatorname{HF}^{\mathrm{bd}}_{j,\Zt}(X;L)\le a_jL+\beta_j$.
Every relative $j$-cycle $z$ that bounds in $(M,\partial M)$
then satisfies
\begin{equation}\label{wz:eq:transfer}
 \inf_{\dr Q=z}\Mass_g(Q)
 \le\Lambda^{2j+1}a_j\Mass_g(z)
       +\frac{\Lambda^{j+1}}2\beta_j.
\end{equation}
\end{proposition}

\begin{proof}
Fix an integral representative $T$ of $z$. By
Lemma~\ref{wz:lem:double-chain}, $\mathcal Dz$ bounds and has mass
at most $2\Lambda^j\Mass_g(z)$. For $\eta>0$ choose
\[
 \partial B=\mathcal Dz,\quad
 \Mass_h(B)\le2\Lambda^ja_j\Mass_g(z)+\beta_j+\eta.
\]
Apply Lemma~\ref{wz:lem:cut} to both halves and pull their
relative fillings back to $(M,g)$. Each has relative
boundary $z$, and their sum of masses is at most
\[
 \Lambda^{j+1}(1+Ct)^{j+1}\Mass_h(B)
 +t(1+Ct)^j\Lambda^{j+1}
   \bigl(\Mass_h((j_+)_\#T)+\Mass_h((j_-)_\#T)\bigr).
\]
The representative $T$ is fixed before taking $t\downarrow0$,
so the second term tends to zero even if the canonical representative
has infinite absolute boundary mass.
Choose the lighter filling and let regular $t\downarrow0$, then
$\eta\downarrow0$, to obtain \eqref{wz:eq:transfer}.
\end{proof}

The same proof gives the nonlinear version
\begin{equation}\label{wz:eq:transfer-general}
 \inf_{\dr Q=z}\Mass_g(Q)\le
 \frac{\Lambda^{j+1}}2
 \operatorname{HF}^{\mathrm{bd}}_{j,\Zt}
       \bigl(X;2\Lambda^j\Mass_g(z)\bigr).
\end{equation}

The factor $1/2$ comes from choosing the lighter filling, not dividing
mod-two chains. No reflection isometry is used; folding the whole
filling onto $M$ would cancel its doubled boundary.

In degree one, completion through a complementary region avoids
choosing a metric on a double.

\begin{lemma}
\label{wz:lem:completion}
Let $(X,h)$ be closed with $H_1(X;\Zt)=0$ and
$\operatorname{HF}^{\mathrm{bd}}_{1,\Zt}(X;L)\le aL+b$.
Let $\Omega\subset X$ be a smooth compact domain such that
$M=X\setminus\Int\Omega$ is connected. Suppose every
$x,y\in\partial\Omega$ can be joined in $\Omega$ by a curve
of length at most $q\,d_M(x,y)$, with an arbitrarily small
additional error allowed. Then
\begin{equation}\label{wz:eq:completion}
 \FH_1(L)\le(1+q)aL+b
\end{equation}
for the induced metric on $M$.
\end{lemma}

\begin{proof}
Fix an integral representative $T$ of $z$. We first pair the points
of $\partial T$ so that the sum of the intrinsic distances in $M$
of the paired points is at most $\Mass(T)$. If $\partial T=0$ the
pairing is empty. Otherwise let $N=\Mass(\partial T)$. By polyhedral
approximation modulo two \cite[Theorem 3.4]{MS2018}, there are
reduced polyhedral chains $P_i$ converging flatly to $T$ in a
shrinking ambient neighborhood of $M$, with
$\Mass(P_i)\le\Mass(T)+o(1)$ and $\Mass(\partial P_i)\le N+o(1)$.
For large $i$, the boundary of $P_i$ has at most $N$ points. Flat
convergence of the boundaries forces a point near each of the $N$
distinct points of $\partial T$, so the endpoints match with
distances tending to zero.

View $P_i$ as a finite graph whose total edge length is $\Mass(P_i)$.
Project each parametrized edge to $M$ by nearest-point projection
to a smooth extension followed by clamping the outer collar to
$\partial M$. On the shrinking neighborhood this map increases the
length of each edge by a factor of at most $1+o(1)$. Join the
projected odd vertices to their matched points of
$\partial T$ by paths in $M$ of total length $o(1)$. In the abstract
graph, pair incident edges at each vertex, leaving one end unpaired
at each prescribed endpoint. The open trails pair the points of
$\partial T$, and their total length is at most
$(1+o(1))\Mass(P_i)+o(1)$, even if projected edges overlap.
Their endpoint distances satisfy the same bound. Since $\partial T$
has only finitely many pairings, one occurs along a subsequence and
has total distance at most $\Mass(T)$.

Join these pairs in $\Omega$ using the hypothesis, with total
additional error at most $\eta>0$. The resulting mod-two chain $C$
satisfies $\partial C=\partial T$ and $\Mass(C)\le q\Mass(T)+\eta$.
Thus $T+C$ is a cycle in $X$ of mass at most $(1+q)\Mass(T)+\eta$
and agrees with $T$ in $\Int M$. It bounds since $H_1(X;\Zt)=0$.
Apply the closed filling estimate to this single cycle and use
Lemma~\ref{wz:lem:cut} on $M$. Letting the filling error and $\eta$
tend to zero gives a relative filling cost at most
$a(1+q)\Mass(T)+b$. Taking the infimum over representatives $T$
and then the supremum over relative cycles of mass at most $L$
gives \eqref{wz:eq:completion}.
\end{proof}

\Needspace{4\baselineskip}
We apply the transfer to the controlled four-dimensional doubles
of \eqref{wz:eq:double-hyp}, without assuming absolute homology vanishing.

\begin{proof}[Proof of Theorem~\ref{thm:relative-wz}]
Relative homology vanishing ensures that every relative
one-cycle bounds. Apply Proposition~\ref{wz:prop:closed}
to $(X,h)$ and Proposition~\ref{wz:prop:transfer} to obtain
\eqref{wz:eq:relative-main}. For $n=4,k=2$,
Theorem~\ref{thm:mass} gives
$\|V\|(M)\le60\FH_1(2D)$, which is
\eqref{wz:eq:mass-main}. The boundary mass assertion follows from
\cite[Corollary~1.5]{PartI}.
\end{proof}

\begin{theorem}
\label{hf:thm:relative}
Let $(M^4,g)$ have a smooth double $X=DM$ with a metric $h$ satisfying
\eqref{wz:eq:double-hyp}. For bounding relative cycles,
\begin{align}
 \FH_2^{\mathrm{bd}}(L)
 &\le\Lambda^5a_2L+\tfrac12\Lambda^3\beta_2,
 \label{hf:eq:relative2}\\
 \FH_3^{\mathrm{bd}}(L)
 &\le\min\{\Vol_g(M)/2,\Lambda^7c_3(d_0)L\}.
 \label{hf:eq:relative3}
\end{align}
For each degree $j=2,3$ this is an estimate for $\FH_j$
provided $H_j(M,\partial M;\Zt)=0$.
No absolute homology hypothesis on $X$ is needed.
\end{theorem}

\begin{proof}
Apply Proposition~\ref{wz:prop:transfer} in degrees two
and three to \eqref{hf:eq:closed2} and the linear part
of \eqref{hf:eq:closed3}. The half-volume bound follows
directly from Proposition~\ref{hf:prop:top} on $M$.
\end{proof}

The homology condition for $\FH_3$ is independent of those for
$\FH_1,\FH_2$: $S^1\times D^3$ has vanishing relative homology
in degrees one and two, but not three.

\begin{remark}\label{wz:rem:scope}
If $\Vol_g(M)\ge v$ and $\diam_g(M)\le d$, the metric comparison
implies
\[
 \Vol_h(X)\ge2\Lambda^{-4}v,\quad
 \diam_h(X)\le2\Lambda d.
\]
Thus these values can be used for $v_0,d_0$ once the Ricci
bound on the smooth double has been established.
Also $2\Lambda^{-1}D\le\diam_h(X)$: any path between the
two copies of a point must reach and leave the seam.
Consequently \eqref{wz:eq:mass-main} is at most
$60\Lambda^4a d_0+30\Lambda^2b$.

Interior Ricci bounds alone do not supply a controlled
smoothing of the double. Theorem~\ref{thm:relative-wz} makes no assertion
that such a smoothing follows from boundary convexity
either. Proposition~\ref{wz:prop:obstruction} shows that the intrinsic
analogue of the filling estimate \eqref{wz:eq:relative-main} fails even
at $L=1$.
\end{remark}

The examples are the following.

\begin{proposition}
\label{wz:prop:obstruction}
There are smooth compact connected four-manifolds $(M_N,g_N)$
with connected smooth boundary such that
\begin{equation}\label{wz:eq:obstruction-geometry}
 \operatorname{Ric}_{g_N}=0,\quad
 \Vol_{g_N}(M_N)=1,\quad \diam_{g_N}(M_N)<15,\quad
 H_1(M_N,\partial M_N;\Zt)=0,
\end{equation}
but $\FH_1^{M_N}(1)\to\infty$ and
$\FH_1^{M_N}(2D_N)\to\infty$, where
$D_N=\sup_{x\in M_N}d_{M_N}(x,\partial M_N)$.
There is also such a sequence with
$H_1(M_N;\mathbb Z)=0$, $|\operatorname{Ric}_{g_N}|\le3$
and $\Vol_{g_N}(M_N)\ge1$, in place of the flatness and
volume-one conditions.
\end{proposition}

\begin{proof}
For $N\ge4$ put $\eps_N=1/N$. Consider the rectangle
$\{\eps_N\le r\le1,\ 0\le\theta\le2\pi N\}$
in the universal cover of the punctured Euclidean plane,
with the flat metric $dr^2+r^2d\theta^2$.
No values of $\theta$ are identified. Round the four
corners inward and smoothly inside coordinate boxes
of radial width $\eps_N/10$ and angular width
$\eps_N/10$, symmetrically under
$\theta\mapsto2\pi N-\theta$. Choose the rounded corners
so that each radial fiber remains an interval containing
$r=2\eps_N$. The resulting surface $\Sigma_N$ is a smooth
disk with flat interior and area
\begin{equation}\label{wz:eq:disk-area}
 A_N=\pi(N-N^{-1})-O(N^{-2})\longrightarrow\infty.
\end{equation}
Its intrinsic diameter is at most $2+4\pi$: join any
two points radially to the circle segment $r=2\eps_N$,
and then along that segment.

The unaffected radial arc $\gamma_N=\{\theta=\pi N\}$ has endpoints
on $\partial\Sigma_N$, length $1-\eps_N$, and separates two regions
of area $A_N/2$. These are its only relative mod-two fillings:
by constancy, two fillings differ by zero or the whole disk. Thus
\begin{equation}\label{wz:eq:disk-fill}
 \operatorname{Fill}_{\Sigma_N}(\gamma_N)=A_N/2.
\end{equation}

Let $T_N$ be a square flat two-torus of area $A_N^{-1}$
and set $M_N=\Sigma_N\times T_N$ with the product metric.
It is flat and has volume one. The rounding can be chosen with
$A_N>1$, and $\diam(T_N)=(2A_N)^{-1/2}$, so
$\diam(M_N)\le\sqrt{(2+4\pi)^2+(2A_N)^{-1}}<15$.
Its boundary is $\partial\Sigma_N\times T_N$,
which is smooth and connected. The manifold retracts
onto $T_N$, so Poincar\'e--Lefschetz duality gives
$H_1(M_N,\partial M_N;\Zt)\simeq H^3(M_N;\Zt)=0$.
Fix $p_N\in T_N$ and let
$z_N=\gamma_N\times\{p_N\}$.
The $1$-Lipschitz projection of pairs $(M_N,\partial M_N)\to
(\Sigma_N,\partial\Sigma_N)$ sends fillings of $z_N$ to fillings
of $\gamma_N$ without increasing mass. Conversely, either half-disk
times $\{p_N\}$ fills $z_N$. Thus \eqref{wz:eq:disk-fill} gives
\[
 \inf_{\dr Q=z_N}\Mass(Q)=A_N/2,\quad
 \Mass(z_N)=1-\eps_N<1.
\]
This proves the claim at argument $1$.

In fact $D_N=(1-\eps_N)/2$. The radial boundary gives
this upper bound everywhere, since rounding only
shortens the radial intervals. At
$(r,\theta)=((1+\eps_N)/2,\pi N)$ the inner and outer
radial boundary arcs have precisely this distance.
A path to either angular side or a rounded corner
has length at least
$\eps_N(\pi N-\eps_N/10)>3$, since $r\ge\eps_N$
along every such path. Hence no such path is shorter.
The product does not change distance to the boundary.
Thus $\Mass(z_N)=2D_N$, proving divergence also at
the argument used in Theorem~\ref{thm:mass}.

For the last assertion replace $T_N$ by the unit
round two-sphere. The product then retracts onto
$S^2$ and has $H_1(M_N;\mathbb Z)=0$ as well as
$H_1(M_N,\partial M_N;\Zt)=0$. Its Ricci eigenvalues
are $0,0,1,1$, its volume is $4\pi A_N\ge1$, and
its diameter is at most
$\sqrt{(2+4\pi)^2+\pi^2}<15$. The same projection
and filling argument applies unchanged.
\end{proof}

Thus two-sided interior Ricci bounds, a positive lower volume bound,
bounded intrinsic diameter and relative homology vanishing do not
give uniform relative filling estimates, even at prescribed total
volume. The closed estimate of Proposition~\ref{wz:prop:closed} is unaffected,
since the boundary geometry is uncontrolled here. The example does not assert divergence of the
least mass of a stationary two-varifold.

\section{Geometric applications}
\label{sec:applications}
\label{sec:width-literature}

In dimension four, Theorem~\ref{thm:relative-wz} bounds the mass of a
two-varifold using $\FH_1$. For hypersurfaces the recursion also uses
$\FH_2$.

\begin{corollary}
\label{hf:cor:hypersurface}
Let $(M^4,g)$ have a smooth double $X=DM$ with a metric $h$ satisfying
\eqref{wz:eq:double-hyp}, and assume
$H_1(M,\partial M;\Zt)=H_2(M,\partial M;\Zt)=0$.
Write $a_1=a$, $\beta_1=b$ for the degree-one constants of
Theorem~\ref{thm:relative-wz}, and let $a_2,\beta_2$ be the degree-two
constants of Proposition~\ref{hf:prop:closed2}.
There is a nonzero integral three-varifold $V$,
stationary with free boundary and almost minimizing
modulo two in small relative annuli, such that
\begin{equation}\label{hf:eq:mass3}
 \begin{split}
 0<\|V\|(M)
 &\le20\FH_2\bigl(3\FH_1(2D)\bigr)\\
 &\le120\Lambda^8a_1a_2D
       +30\Lambda^7a_2\beta_1+10\Lambda^3\beta_2.
 \end{split}
\end{equation}
If $H_{\partial M}$ is nowhere zero, then
$\|V\|(\partial M)=0$. If $H_{\partial M}>0$, its support
is a finite union of smooth compact properly embedded
free-boundary minimal hypersurfaces, with positive integer
multiplicities. Their total area counted with these
multiplicities obeys \eqref{hf:eq:mass3}.
\end{corollary}

\begin{proof}
Theorem~\ref{thm:mass} with $n=4,k=3$ gives the first
inequality, since $(n+1)!/k!=20$. The two homology assumptions
allow us to use both filling estimates for every cycle. More explicitly,
\[
 \begin{aligned}
 \FH_1(2D)&\le2\Lambda^3a_1D+\tfrac12\Lambda^2\beta_1,\\
 \FH_2\bigl(3\FH_1(2D)\bigr)
 &\le3\Lambda^5a_2\bigl(2\Lambda^3a_1D+
       \tfrac12\Lambda^2\beta_1\bigr)+\tfrac12\Lambda^3\beta_2.
 \end{aligned}
\]
Multiplication by $20$ gives \eqref{hf:eq:mass3}.
Since $\operatorname{tr}_{T\partial M}\SFF_{\partial M}=H_{\partial M}$,
the boundary mass and proper embeddedness assertions follow from
\cite[Corollaries~1.5 and~1.6]{PartI}.
\end{proof}

These hypersurfaces may have empty boundary. Hypersurface regularity
does not apply to the two-varifolds of Theorem~\ref{thm:relative-wz}
in a four-manifold, and no smoothness is asserted there. The separate
top-degree estimate for $\FH_3$ is not required for $k=3$.

\Needspace{4\baselineskip}
Product collars provide a direct source of controlled doubles.

\begin{corollary}\label{wz:cor:product}
Suppose $(M^4,g)$ has a product metric
$g=ds^2+g_{\partial M}$ in a boundary collar and satisfies
\[
 |\operatorname{Ric}_g|\le3,\quad
 \Vol_g(M)\ge v,\quad \diam_g(M)\le d,\quad
 H_1(M,\partial M;\Zt)=0.
\]
Then
\begin{align}
 \FH_1(L)&\le a(2v,2d)L+\tfrac12b(2v,2d),
 \label{wz:eq:product-fill}\\
 0<\|V\|(M)&\le120a(2v,2d)D+30b(2v,2d)
 \label{wz:eq:product-mass}
\end{align}
for a varifold with the stationarity and almost-minimizing
properties of Theorem~\ref{thm:relative-wz}.
\end{corollary}

\begin{proof}
Reflection gives $DM$ a smooth metric with the same Ricci bound,
because the metric is independent of the normal coordinate near the
seam. Its volume is twice that of $M$. Points in the same half can
be joined at cost at most $d$; points in opposite halves can be
joined through any fixed seam point at cost at most $2d$.
Apply Theorem~\ref{thm:relative-wz} with $\Lambda=1$;
no absolute homology assumption is needed.
\end{proof}

The relative homology condition here is strictly weaker
than $H_1(DM;\mathbb Z)=0$. For example, the topology
$M=S^1\times B^3$ has $H_1(M,\partial M;\Zt)=0$, while
$DM=S^1\times S^3$ has nonzero first homology. Product
collars have totally geodesic boundary, so
Corollary~\ref{wz:cor:product} alone does not imply zero
boundary mass.

For punctured closed manifolds, completion through the removed ball
avoids smoothing a double. Here the boundary has strictly negative
two-plane traces.

\begin{corollary}
\label{wz:cor:puncture}
Let $(X^4,h)\in\mathcal C(v_0,d_0)$ and
$H_1(X;\Zt)=0$. There is $r_*(X,h)>0$ such that for every
$p\in X$ and $0<r<r_*$ the domain
$M_r=X\setminus B_r(p)$ satisfies
\begin{equation}\label{wz:eq:puncture-fill}
 \FH_1^{M_r}(L)\le2a(v_0,d_0)L+b(v_0,d_0).
\end{equation}
It admits a nonzero integral two-varifold, stationary
with free boundary and almost minimizing modulo two
in small relative annuli, such that
\begin{equation}\label{wz:eq:puncture-mass}
 0<\|V\|(\Int M_r)=\|V\|(M_r)
 \le240a(v_0,d_0)D_r+60b(v_0,d_0),
\end{equation}
where
\begin{equation}\label{wz:eq:puncture-distance}
 D_r=\sup_{x\in M_r}d_{M_r}(x,\partial M_r)
 =\max_{x\in X}d_X(x,p)-r\le d_0-r.
\end{equation}
In particular, the right-hand side can be bounded
independently of $p,r$ by $240a(v_0,d_0)d_0+60b(v_0,d_0)$.
\end{corollary}

\begin{proof}
Compactness gives $r_*(X,h)>0$, uniform in $p$, so each
$\overline{B_r(p)}$, $r<r_*$, is a strongly convex smooth ball with
strictly convex boundary. No uniformity over $\mathcal C(v_0,d_0)$
is claimed. In the exact sequence of $(X,\overline{B_r(p)})$, the
terms $H_1(X;\Zt)$ and the kernel of
$H_0(\overline{B_r(p)};\Zt)\to H_0(X;\Zt)$ both vanish.
Thus $H_1(X,\overline{B_r(p)};\Zt)=0$; excision through a collar
identifies this group with $H_1(M_r,\partial M_r;\Zt)$.
Boundary points can be joined inside
the ball by an ambient minimizing geodesic, no longer than their
intrinsic distance in $M_r$. Apply Lemma~\ref{wz:lem:completion} with
$q=1$, Proposition~\ref{wz:prop:closed}, and Theorem~\ref{thm:mass}.

A minimizing geodesic from $x$ to $p$, stopped at the sphere, has
length $d_X(x,p)-r$ and lies in $M_r$; the triangle inequality gives
the reverse bound. This proves \eqref{wz:eq:puncture-distance}.
The outward normal points into the removed ball, so
$\SFF_{\partial M_r}$ is negative definite. The zero boundary mass
assertion follows from \cite[Corollary~1.5]{PartI}.
\end{proof}

These bounds assert neither smoothness nor nonempty boundary. The
puncture hypotheses do not supply the controlled double required in
Corollary~\ref{hf:cor:hypersurface}.

\begin{remark}[Comparison with width bounds]
\label{width:rem:comparison}
For hypersurfaces there is a different route to mass bounds, which uses
the width of a single sweepout instead of filling functions. Let
$(X^n,h)$ be a closed connected manifold with
$\operatorname{Ric}_h\ge-(n-1)a^2$, $a\ge0$. Glynn-Adey and Liokumovich
\cite[Corollary 1.4]{GAL2017} bound its first mod-two hypersurface
width by
\[
 W(X,h)\le C_n\max\{1,a\Vol_h(X)^{1/n}\}\Vol_h(X)^{(n-1)/n},
\]
and \cite[Theorem 5.1]{GAL2017} gives a corresponding estimate for a
smooth compact domain $U\subset X$ in terms of $\Vol(U)$ and
$\operatorname{Area}(\partial U)$; for the restriction of sweepouts to
domains, see also \cite[Lemma 2.15]{LMN2018}. Relative min--max theory
\cite{LZ2021,GLWZ2021} converts a sweepout of $U$ by relative cycles
into a nonzero free-boundary stationary integral $(n-1)$-varifold in
$U$, almost minimizing modulo two in small relative annuli, whose mass
is at most the maximal mass of the sweepout; for $3\le n\le7$ it is
supported on smooth almost properly embedded free-boundary minimal
hypersurfaces. No vanishing of $H_j(U,\partial U;\Zt)$ is needed.

These varifolds also satisfy the boundary conclusions of
\cite[Corollaries~1.5 and~1.6]{PartI}. Under its homology hypotheses,
Theorem~\ref{thm:mass} gives mass bounds for $k$-varifolds
for every $2\le k\le n-1$, in terms of $D$ and the relative filling
functions. In dimension four, the estimates of
Sections~\ref{sec:closed-fill}--\ref{sec:relative-wz} additionally
control the filling mass of every bounding cycle in the relevant
degrees. A finite width does not imply a finite all-cycle filling
function: $S^2\times D^2\subset S^2\times S^2$ has finite relative
width, but $H_2(S^2\times D^2,\partial(S^2\times D^2);\Zt)\ne0$, so its
function $\FH_2$ is infinite at sufficiently large arguments. None of
these existence arguments ensures that the resulting hypersurface has
nonempty boundary.
\end{remark}

\appendix
\section{Proof of Theorem~\ref{thm:relative-selection}\textup{(ii)}}
\label{app:selection}

We use the notation of Theorem~\ref{thm:relative-selection} and prove
part \textup{(ii)} with $N(m)=2\cdot3^m-1$. We follow the min--max
argument of Pitts \cite[Theorem~4.10]{Pitts1981} and its adaptation
to free-boundary hypersurfaces by Li and Zhou \cite[Theorem~4.21]{LZ2021}.
Proposition~\ref{app:tightening} discretizes the pull-tight deformation
for relative mod-two cycles of dimension $k$, fixing the prescribed
parameter faces and keeping the mass below $\Lambda$.
Lemma~\ref{app:selection-lemma} selects disjoint annuli, and
Lemma~\ref{app:combination} combines the deformations supported in them
to obtain $\mathcal P_{N(m)}$. For the corresponding arguments for
closed cycles and for isotopies of surfaces, see
\cite[Proposition~4.9 and pp.~165--174]{Pitts1981} and
\cite[Appendix]{CGK2018}.

\subsection{Discrete families}

For $j$ at least the subdivision level of $X$, recall that $X(j)_0$ is the set of vertices of the subdivision $X(j)$ of $X$ into cubes of side $3^{-j}$, and put $A(j)_0=A\cap X(j)_0$. Two vertices are \emph{adjacent} if they are the endpoints of an edge of $X(j)$. A \emph{discrete homotopy} is a map
\[
 \psi:\{0,1,\ldots,P\}\times X(j)_0\longrightarrow\ZC_k(M,\partial M;\Zt).
\]
Its fineness $\mathbf f(\psi)$ is the largest value of $\Mass(\psi(i,v)-\psi(i',v'))$ over all pairs such that either $i=i'$ and $v,v'$ are adjacent, or $v=v'$ and $|i-i'|=1$. We say that $\psi$ vanishes on $A$ if $\psi(i,v)=0$ for all $i$ and all $v\in A(j)_0$.

For $j'>j$, let $n:X(j')_0\to X(j)_0$ assign to $w$ the point of the grid of mesh $3^{-j}$ that is nearest to $w$ in the supremum norm $|\cdot|_\infty$ of $\R^m$. This point is unique because $3^{j'-j}$ is odd.

\begin{lemma}\label{app:nearest}
Let $j'>j$.
\begin{enumerate}[label=\textup{(\roman*)},leftmargin=2.2em]
\item If $w\in X(j')_0$ lies in a cell $\alpha$ of $X(j)$, then $n(w)$ is a vertex of $\alpha$. In particular, $n$ maps $X(j')_0$ into $X(j)_0$ and $A(j')_0$ into $A(j)_0$.
\item If $w,w'$ are adjacent in $X(j')$, then $n(w)$ and $n(w')$ are equal or adjacent in $X(j)$.
\item Let $\psi$ be a discrete homotopy on $X(j)_0$ with $P\le3^{j'}$. Then
\[
 \psi'(t,w)=\psi\bigl(\min\{3^{j'}t,P\},n(w)\bigr)
\]
is a map on the vertex set of $(I\times X)(j')$ with $\mathbf f(\psi')\le\mathbf f(\psi)$, and it has the same set of values as $\psi$.
\end{enumerate}
\end{lemma}

\begin{proof}
For \textup{(i)}, each coordinate of $w$ lies in the corresponding coordinate interval of $\alpha$, which is a grid value or an interval between two consecutive grid values; the nearest grid value is therefore an endpoint of this interval. For \textup{(ii)}, the edge joining $w$ and $w'$ lies in a cell $\alpha$ of $X(j)$, and $w,w'$ differ in a single coordinate. By \textup{(i)}, $n(w)$ and $n(w')$ are vertices of $\alpha$ that differ at most in this coordinate, so they are equal or joined by an edge of $\alpha$. Part \textup{(iii)} follows from \textup{(ii)}, since adjacent times correspond to indices that differ by at most one.
\end{proof}

The next lemma returns from discrete homotopies to the class $\mathcal C$.

\begin{lemma}\label{app:return}
Let $\delta_0,C_0$ be the constants of Theorem~\ref{inp:interpolation}\textup{(ii)} for $I\times X$, and let $a_X,C_X$ be those of Lemma~\ref{mass:small-cone} for $(X,A)$. Let $\Phi\in\mathcal C$, and let $\theta$ be a map from the vertex set of $(I\times X)(j)$ to $\ZC_k(M,\partial M;\Zt)$ which vanishes at the vertices in $I\times A$ and has $\mathbf f(\theta)<\delta_0$. Suppose that, for some $\omega\ge0$,
\[
\begin{aligned}
 \Mass(\theta(0,v)-\Phi(v))&\le\omega&&\text{for every }v\in X(j)_0,\\
 \Mass(\Phi(x)-\Phi(y))&\le\omega&&\text{for all $x,y$ in a common cell of }X(j),
\end{aligned}
\]
and that $e=C_0\mathbf f(\theta)+2\omega$ satisfies
\[
 e<a_X,\quad
 \max\Mass(\theta)+C_0\mathbf f(\theta)<\Lambda,\quad
 \sup_X\Mass(\Phi)+(1+C_X)e<\Lambda.
\]
Then the Almgren extension $\Theta$ of $\theta$ satisfies $\Theta(1,\cdot)\in\mathcal C$ and
\[
 \sup_{x\in X}\Mass(\Theta(1,x))\le\max_{v\in X(j)_0}\Mass(\theta(1,v))+C_0\mathbf f(\theta).
\]
\end{lemma}

\begin{proof}
The map $\Theta$ is mass-continuous on $I\times X$ and agrees with $\theta$ at the vertices. As in the proof of Lemma~\ref{mass:small-cone}, the extension is chosen to be zero on every cell all of whose vertex values vanish, so $\Theta$ vanishes on $I\times A$. By \eqref{inp:eq:almgren-extension}, every value of $\Theta$ differs in mass by at most $C_0\mathbf f(\theta)$ from the value at a vertex of the same cell. This gives the last inequality and also $\sup\Mass(\Theta)<\Lambda$.

Put $B=\Theta(0,\cdot)$. For $x$ in a cell of $X(j)$ with vertex $v$,
\[
 \Mass(B(x)-\Phi(x))\le\Mass(B(x)-B(v))+\Mass(\theta(0,v)-\Phi(v))+\Mass(\Phi(v)-\Phi(x))\le e.
\]
The map $\Phi-B$ vanishes on $A$ and has mass at most $e<a_X$, so Lemma~\ref{mass:small-cone} gives a mass-continuous homotopy $G$ from $\Phi-B$ to zero, vanishing on $A\times I$, with $\sup\Mass(G)\le C_Xe$. Then $(x,t)\mapsto B(x)+G(x,t)$ joins $\Phi$ to $B$ relative to $A$, and its masses are at most $\sup_X\Mass(\Phi)+(1+C_X)e<\Lambda$. Following a homotopy from $\Phi_0$ to $\Phi$ as in the definition of $\mathcal C$ by this homotopy and then by $t\mapsto\Theta(t,\cdot)$, we see that $\Theta(1,\cdot)\in\mathcal C$.
\end{proof}

\subsection{Tightening}

\begin{proposition}\label{app:tightening}
Let $\Phi_i\in\mathcal C$ satisfy $\sup_X\Mass(\Phi_i)\to a$. After passing to a subsequence, there are integers $j_i\to\infty$ and discrete homotopies
\[
 \psi_i:\{0,\ldots,P_i\}\times X(j_i)_0\longrightarrow\ZC_k(M,\partial M;\Zt),
\]
vanishing on $A$, with the following properties.
\begin{enumerate}[label=\textup{(\roman*)},leftmargin=2.2em]
\item $\psi_i(0,v)=\Phi_i(v)$ for every $v\in X(j_i)_0$, and $\Mass(\Phi_i(x)-\Phi_i(y))\le1/i$ whenever $x,y$ lie in a common cell of $X(j_i)$.
\item $\mathbf f(\psi_i)\to0$ and $\max\Mass(\psi_i)\le\sup_X\Mass(\Phi_i)+o(1)$.
\item Put $\phi_i=\psi_i(P_i,\cdot)$. If $i_\ell\to\infty$, $v_\ell\in X(j_{i_\ell})_0$, $\Mass(\phi_{i_\ell}(v_\ell))\to a$ and $|\phi_{i_\ell}(v_\ell)|\rightharpoonup V$, then $V$ is stationary with free boundary.
\end{enumerate}
\end{proposition}

\begin{proof}
On the compact space of $k$-varifolds of mass at most $\Lambda$, let
$D_{t,W}$, $0\le t\le1$, denote the diffeomorphisms of $M$ given by
the pull-tight construction \cite[Theorem~4.3]{Pitts1981},
\cite[proof of Proposition~4.17]{LZ2021}. They depend continuously on
$(t,W)$ in $C^1$, start at the identity, preserve $\partial M$, and
do not increase the mass of $W$ along the flow. For stationary $W$
they are the identity. The mass drop
\[
 b(W)=\|W\|(M)-\|(D_{1,W})_\#W\|(M)
\]
is continuous and strictly positive for every $W$ that is not
stationary with free boundary. Near each nonstationary $W$ one chooses a field
tangent to $\partial M$ with negative first variation, combines these
fields by a partition of unity on varifold space, and chooses a
continuous flow time that tends to zero at the stationary set. The
construction applies to $k$-planes in every dimension under consideration.

Fix $i$, write $Y=I\times X$ and $B=(\{0\}\times X)\cup(I\times A)$,
and consider the relative cycles
\[
 F_i(t,x)=(D_{t,|\Phi_i(x)|})_\#\Phi_i(x).
\]
The diffeomorphisms preserve the boundary, so they carry canonical
representatives to canonical representatives. Thus $F_i$ is continuous
both in the relative flat distance and as a family of varifolds, and
$\Mass(F_i(t,x))\le\Mass(\Phi_i(x))$. Its restriction to $B$ is
mass-continuous. The family has no concentration of mass since its weight
measures form a weakly compact family and each is atomless for
$k\ge1$. Otherwise, passing to convergent subsequences of the
parameters and the centers of shrinking balls would give an atom in a
limiting weight measure.

We use the cellwise discretization in \cite[Appendix~A]{LMN2018}.
For this fixed $i$ it gives coarse levels $h_\ell\to\infty$, finer
levels $j_\ell\ge h_\ell$, errors $e_\ell\downarrow0$, and maps
$\theta_\ell:Y(j_\ell)_0\to\ZC_k(M,\partial M;\Zt)$ satisfying
\[
 \mathbf f(\theta_\ell)+
 \sup_{z\in Y(j_\ell)_0}\Flat(\theta_\ell(z)-F_i(z))\le e_\ell,
 \quad \theta_\ell=F_i\quad\text{on }B(j_\ell)_0.
\]
For each \emph{coarse cell} $\alpha$ of $Y(h_\ell)$, the mass satisfies
\begin{equation}\label{app:eq:tightening-cell}
 \max_{z\in\alpha\cap Y(j_\ell)_0}\Mass(\theta_\ell(z))
       \le \sup_{y\in\alpha}\Mass(F_i(y))+e_\ell .
\end{equation}
To see this, use the induction over the parameter skeleton in
\cite[Section~13.7, condition~(2) and (79)]{MN2014}.
For relative mod-two $k$-cycles the required extension on a single
cell is \cite[Lemma~A.1]{LMN2018}. At tolerance $\delta$, each extension
increases the maximum of the already assigned boundary masses by at
most $\delta$. Induction bounds the mass on a coarse $p$-cell by the
maximum of its corner masses plus $p\delta$. Shared faces are assigned
once and pulled back by the nearest-vertex maps at subsequent stages.
On a cell of $B$ use the coarse samples pulled back by the same maps
throughout the induction. These satisfy the same fineness and local
flat estimates by choosing the coarse grid with mass oscillation less
than $\delta$ on each such cell. At the end replacing the pulled-back
values on $B$ by the exact samples changes each vertex by at most
$\delta$ in mass and hence each edge bound by at most $2\delta$.
The fixed-face construction in \cite[proof of Theorem~13.1]{MN2014}
applies to $B$ because every face of a cell in $B$ also lies in $B$
and the prescriptions agree on $\{0\}\times A$. Letting $\delta$ tend to zero
gives the estimates above.

For fixed $i$, continuity of $y\mapsto\Mass(F_i(y))$ makes the
oscillation of this function on the coarse cells tend uniformly to
zero. We claim that
\[
 \sup_{z\in Y(j_\ell)_0}\dV\bigl(|\theta_\ell(z)|,|F_i(z)|\bigr)
 \longrightarrow0.
\]
Indeed, along any sequence $z_\ell\to z$ the cycles converge flatly
to $F_i(z)$, and their masses converge to $\Mass(F_i(z))$ by
\eqref{app:eq:tightening-cell} and lower semicontinuity. To identify
any varifold limit $U$, put $T=F_i(z)$ and vary the metric to
$g_u=g+u h$ for a smooth symmetric tensor $h$ and $|u|$ small.
The relative flat topologies agree. Lower semicontinuity for $g_u$
gives $\int J_{g_u}(P)\,dU\ge\Mass_{g_u}(T)$, where $J_{g_u}(P)$
is the ratio of the $g_u$- and $g$-volume elements on $P$.
Equality holds at $u=0$, so differentiation gives
$\int\operatorname{tr}_P h\,dU=\int\operatorname{tr}_P h\,d|T|$.
Taking $h=\varphi g$ gives equality of the weights. For almost every
$x$ with respect to this weight, testing all $h$ shows that the
average of the orthogonal projections onto the planes of $U$ at $x$
is the projection onto the approximate tangent plane of $T$ at $x$.
The trace of this average on the perpendicular complement is zero.
Each projection has nonnegative trace there, so almost every plane
equals the tangent plane. Thus $U=|T|$. Compactness of $Y$ gives
uniformity, as in \cite[proof of Theorem~4.5, (4.1)--(4.2)]{GLWZ2021}.

Now choose a refinement separately for each $i$ so that the fineness,
the mass error in \eqref{app:eq:tightening-cell}, and the uniform
varifold error are at most $1/i$. Require also $j_i\ge i$ and mass
oscillation at most $1/i$ for $\Phi_i$ on cells of $X(j_i)$.
Reading the time layers of $\theta_\ell$ as $\psi_i$ gives
\textup{(i)} and \textup{(ii)}. For a subsequence as in
\textup{(iii)}, put $W_i=|\Phi_i(v_i)|$ and pass to a further
subsequence with $W_i\rightharpoonup W$. Uniform varifold approximation
on the last row gives $\|(D_{1,W_i})_\#W_i\|(M)\to a$. Since
$\|W_i\|(M)\le\sup_X\Mass(\Phi_i)\to a$ and $b(W_i)\ge0$, we have
$b(W)=0$. Thus $W$ is stationary and $D_{1,W}$ is the identity.
Continuity of pushforward gives $(D_{1,W_i})_\#W_i\rightharpoonup W$,
so the prescribed limit $V$ equals $W$.
\end{proof}

\subsection{Selection of annuli}

Recall that $\An(p;\sigma,\rho)=\{x\in M:\sigma<d(x,p)<\rho\}$, where $d$ is the distance in the fixed smooth extension of $(M,g)$.

\begin{lemma}\label{app:selection-lemma}
Let $G$ be a finite graph in which every vertex has at most $K$ neighbors, and let $N\ge2K+1$. Suppose that every vertex $v$ of $G$ is given a point $p_v\in M$ and numbers $0<\sigma^v_j<\rho^v_j$, $1\le j\le N$, such that
\begin{equation}\label{app:eq:separation}
 3\rho^v_{j+1}<\sigma^v_j\quad(1\le j<N).
\end{equation}
Then one can choose indices $j(v)$ such that
\[
 \An\bigl(p_v;\sigma^v_{j(v)},\rho^v_{j(v)}\bigr)\cap
 \An\bigl(p_{v'};\sigma^{v'}_{j(v')},\rho^{v'}_{j(v')}\bigr)=\varnothing
 \quad\text{whenever $v$ and $v'$ are neighbors.}
\]
\end{lemma}

\begin{proof}
For neighbors $v,v'$ put $D=d(p_v,p_{v'})$, and say that the index $j$ of $v$ is \emph{blocked} by $v'$ if $\sigma^v_j/2<D<3\rho^v_j/2$. By \eqref{app:eq:separation}, the intervals $(\sigma^v_j/2,3\rho^v_j/2)$, $1\le j\le N$, are pairwise disjoint, so each neighbor blocks at most one index of $v$. Call an index of $v$ \emph{admissible} if no neighbor blocks it. Every vertex has at least $N-K\ge K+1$ admissible indices.

We first prove the following claim: if $j$ is admissible for $v$, if $j'$ is admissible for a neighbor $v'$ of $v$, and if $\rho^{v'}_{j'}\le\sigma^v_j/2$, then the annuli $\An(p_v;\sigma^v_j,\rho^v_j)$ and $\An(p_{v'};\sigma^{v'}_{j'},\rho^{v'}_{j'})$ are disjoint. Write $\sigma=\sigma^v_j$, $\rho=\rho^v_j$ and $\rho'=\rho^{v'}_{j'}$, and suppose that $x$ lies in both annuli. Then
\[
 \sigma<d(x,p_v)\le d(x,p_{v'})+D<\rho'+D,
 \quad
 D\le d(p_v,x)+d(x,p_{v'})<\rho+\rho'.
\]
Since $\rho'\le\sigma/2<\rho/2$, this gives $\sigma/2<D<3\rho/2$, so $j$ is blocked by $v'$, contrary to the assumption.

We now choose the indices one vertex at a time. At each stage, every vertex that has not yet been treated carries a set of \emph{available} indices; initially these are its admissible indices. A stage consists of the following. Among all pairs $(v,j)$ such that $v$ is untreated and $j$ is available for $v$, pick one for which $\rho^v_j$ is smallest, set $j(v)=j$, and declare $v$ treated. Then, for every untreated neighbor $v'$ of $v$, remove from the available indices of $v'$ all $j'$ with $\sigma^{v'}_{j'}/2<\rho^v_j$. By the choice of $(v,j)$, every index $j'$ that is available for $v'$ at this stage has $\rho^{v'}_{j'}\ge\rho^v_j$, so a removed index satisfies
$\rho^v_j\in(\sigma^{v'}_{j'}/2,\rho^{v'}_{j'}]$.
By \eqref{app:eq:separation} these intervals are pairwise disjoint for different $j'$, and therefore at most one index of $v'$ is removed. Consequently a vertex loses at most one available index for each of its neighbors, hence at most $K$ in all, and it always keeps at least one available index. The procedure therefore treats every vertex.

Finally, let $v,v'$ be neighbors, and suppose that $v$ is treated before $v'$. The index $j(v')$ was not removed at the stage of $v$, so $\rho^v_{j(v)}\le\sigma^{v'}_{j(v')}/2$. Both indices are admissible, and the claim, applied with the roles of $v$ and $v'$ interchanged, shows that the two annuli are disjoint.
\end{proof}

\subsection{Combination of local deformations}

\begin{lemma}\label{app:combination}
Let $\phi:X(j)_0\to\ZC_k(M,\partial M;\Zt)$ vanish on $A(j)_0$, and put $\delta=\mathbf f(\phi)$. Let $W\subset X(j)_0$ be a set of vertices with $\Mass(\phi(v))>m\delta$ for all $v\in W$, and let $\delta_1,\zeta>0$ and $q\ge1$. Suppose that every $v\in W$ is given a relatively open set $U_v\subset M$ and relative cycles $\phi(v)=T^v_0,T^v_1,\ldots,T^v_q$ such that
\begin{equation}\label{app:eq:local-sequence}
\begin{gathered}
 \spt(T^v_l-\phi(v))\subset U_v,\quad
 \Mass(T^v_l-T^v_{l-1})\le\delta_1,\\
 \Mass(T^v_l)\le\Mass(\phi(v))+\delta_1,\quad
 \Mass(T^v_q)\le\Mass(\phi(v))-\zeta,
\end{gathered}
\end{equation}
and that $U_v\cap U_{v'}=\varnothing$ whenever $v,v'\in W$ are distinct and $|v-v'|_\infty\le3^{-j}$. Let $j'>j$ satisfy $3^{j'-j}\ge2q$, let $n:X(j')_0\to X(j)_0$ be the nearest-vertex map, and put $\kappa=2^m(2m\delta+\delta_1)$. Then there is a discrete homotopy
$h:\{0,\ldots,q\}\times X(j')_0\to\ZC_k(M,\partial M;\Zt)$,
vanishing on $A$, such that for all $i$ and $w$:
\begin{enumerate}[label=\textup{(\roman*)},leftmargin=2.2em]
\item $h(0,w)=\phi(n(w))$;
\item $\mathbf f(h)\le\delta+2^{m+1}\delta_1$;
\item $\Mass(h(i,w))\le\Mass(\phi(n(w)))+\kappa$;
\item $\Mass(h(q,w))\le\Mass(\phi(n(w)))-\zeta+\kappa$ if $n(w)\in W$.
\end{enumerate}
\end{lemma}

The idea is to add to the cycle at $w$ the changes $T^v_l-T^v_0$ of the nearby vertices $v\in W$. The full change of $v$ is applied on the block of vertices $w$ with $n(w)=v$, and the change is switched off gradually, one term of the sequence at a time, over the adjacent blocks. Adding the same change to two cycles does not alter their difference, which keeps the fineness small; and changes that are active at the same vertex $w$ have disjoint supports, which keeps the mass under control.

\begin{proof}
For $v\in X(j)_0$, let $\operatorname{St}(v)$ be the union of the cells of $X(j)$ that contain $v$. Put $b=(3^{j'-j}-1)/2$; by Lemma~\ref{app:nearest}\textup{(i)}, $n(w)=v$ exactly when $w\in\operatorname{St}(v)$ and $3^{j'}|w-v|_\infty\le b$. For $v\in W$ define $\lambda_v:X(j')_0\to\{0,\ldots,q\}$ by
\[
 \lambda_v(w)=\max\bigl\{0,\ q-\max\{0,\,3^{j'}|w-v|_\infty-b\}\bigr\}
 \quad\text{for }w\in\operatorname{St}(v),
\]
and $\lambda_v(w)=0$ for $w\notin\operatorname{St}(v)$. Thus $\lambda_v(w)=q$ if $n(w)=v$.

We claim that $|\lambda_v(w)-\lambda_v(w')|\le1$ for adjacent $w,w'$. If both lie in $\operatorname{St}(v)$, this holds because $3^{j'}|w-v|_\infty$ and $3^{j'}|w'-v|_\infty$ differ by at most one. Suppose that $w\in\operatorname{St}(v)$ and $w'\notin\operatorname{St}(v)$. The edge joining them lies in a cell $\beta$ of $X(j)$ that does not contain $v$, and $w$ lies in a cell $\alpha$ that contains $v$. Hence $w$ lies in $\alpha\cap\beta$, which is a face of $\alpha$ not containing $v$. Every point of such a face differs from $v$ by exactly $3^{-j}$ in some coordinate, so $|w-v|_\infty=3^{-j}$. Therefore $3^{j'}|w-v|_\infty-b=(3^{j'-j}+1)/2>q$ and $\lambda_v(w)=0=\lambda_v(w')$.

For $w\in X(j')_0$ let $\operatorname{Act}(w)=\{v\in W:\lambda_v(w)>0\}$. We record two properties of this set.
\begin{enumerate}[label=\textup{(\alph*)},leftmargin=2.2em]
\item If $v\in\operatorname{Act}(w)$, then $3^{j'}|w-v|_\infty<q+b$, and $q+b<3^{j'-j}$ because $2q\le3^{j'-j}$. Thus $|w-v|_\infty<3^{-j}$: every coordinate of $v$ is one of the at most two values of the grid of mesh $3^{-j}$ at distance less than $3^{-j}$ from the corresponding coordinate of $w$. Hence $\operatorname{Act}(w)$ has at most $2^m$ elements, and any two of them are at distance at most $3^{-j}$ in the supremum norm, so that the sets $U_v$, $v\in\operatorname{Act}(w)$, are pairwise disjoint.
\item If $v\in\operatorname{Act}(w)$, then $w$ lies in a cell of $X(j)$ that contains $v$, and $n(w)$ is a vertex of this cell by Lemma~\ref{app:nearest}\textup{(i)}. Two vertices of a cell are joined by at most $m$ of its edges, so $\Mass(\phi(n(w))-\phi(v))\le m\delta$.
\end{enumerate}

Put $\Delta^v_l=T^v_l-T^v_0$, so that $\Delta^v_0=0$, and define
\[
 h(i,w)=\phi(n(w))+\sum_{v\in W}\Delta^v_{\min\{i,\lambda_v(w)\}}.
\]
Only the vertices $v\in\operatorname{Act}(w)$ contribute to the sum. Property \textup{(i)} is clear. If $w\in A(j')_0$, then $n(w)\in A(j)_0$ and $\phi(n(w))=0$; moreover $\operatorname{Act}(w)$ is empty, since $v\in\operatorname{Act}(w)$ would give $\Mass(\phi(v))\le m\delta$ by \textup{(b)}. Thus $h$ vanishes on $A$.

For the fineness, consider first $h(i+1,w)-h(i,w)$. Each $v\in\operatorname{Act}(w)$ contributes zero or a difference $T^v_{l+1}-T^v_l$, so the mass is at most $2^m\delta_1$ by \textup{(a)}. Next let $w,w'$ be adjacent. Then
$h(i,w)-h(i,w')$ is the sum of $\phi(n(w))-\phi(n(w'))$, which has mass at most $\delta$ by Lemma~\ref{app:nearest}\textup{(ii)}, and of the terms $\Delta^v_{\min\{i,\lambda_v(w)\}}-\Delta^v_{\min\{i,\lambda_v(w')\}}$ with $v\in\operatorname{Act}(w)\cup\operatorname{Act}(w')$. By the claim, the two indices in each term differ by at most one, so each term has mass at most $\delta_1$, and by \textup{(a)} there are at most $2\cdot2^m$ terms. This proves \textup{(ii)}.

For the masses, fix $(i,w)$ and put $u=n(w)$, $S=\phi(u)$ and $l_v=\min\{i,\lambda_v(w)\}$. The cycles $\Delta^v_{l_v}$, $v\in\operatorname{Act}(w)$, have compact supports in the pairwise disjoint open sets $U_v$. Hence $h(i,w)$ coincides with $S$ near every point outside these sets and with $S+\Delta^v_{l_v}$ in $U_v$, and therefore
\begin{equation}\label{app:eq:mass-splitting}
 \Mass(h(i,w))=\|S\|\Bigl(M\setminus\bigcup_{v}U_v\Bigr)
 +\sum_{v\in\operatorname{Act}(w)}\|S+\Delta^v_{l_v}\|(U_v).
\end{equation}
Fix $v\in\operatorname{Act}(w)$ and write $l=l_v$. Since $T^v_l$ coincides with $\phi(v)$ outside a compact subset of $U_v$,
\[
 \|T^v_l\|(U_v)=\Mass(T^v_l)-\|\phi(v)\|(M\setminus U_v)
 \le\|\phi(v)\|(U_v)+\delta_1 .
\]
Using $S+\Delta^v_l=T^v_l+(S-\phi(v))$ and \textup{(b)}, we obtain
\[
 \|S+\Delta^v_l\|(U_v)\le\|\phi(v)\|(U_v)+\delta_1+m\delta
 \le\|S\|(U_v)+2m\delta+\delta_1 .
\]
Inserting this into \eqref{app:eq:mass-splitting} and using that $\operatorname{Act}(w)$ has at most $2^m$ elements gives \textup{(iii)}. Finally, let $u\in W$ and $i=q$. Then $u\in\operatorname{Act}(w)$, $l_u=q$ and $S+\Delta^u_q=T^u_q$, and the last inequality in \eqref{app:eq:local-sequence} gives
$\|T^u_q\|(U_u)=\Mass(T^u_q)-\|S\|(M\setminus U_u)\le\|S\|(U_u)-\zeta$.
Using this estimate for $v=u$ and the previous one for the other elements of $\operatorname{Act}(w)$ proves \textup{(iv)}.
\end{proof}

\subsection{Conclusion of the proof}

\begin{proof}[Proof of Theorem~\ref{thm:relative-selection}\textup{(ii)}]
\noindent\textbf{Step 1: A critical sequence.}\enspace
Choose $\Phi_i\in\mathcal C$ with $\sup_X\Mass(\Phi_i)\to a$, let $\psi_i$ and $\phi_i=\psi_i(P_i,\cdot)$ be given by Proposition~\ref{app:tightening}, and put $L_i=\max_v\Mass(\phi_i(v))$. By part \textup{(ii)} of the proposition, $L_i\le a+o(1)$. For the reverse inequality, regard $\psi_i$ as a map $\theta_i$ on the vertex set of $(I\times X)(j)$ for a large $j$, as in Lemma~\ref{app:nearest}\textup{(iii)}. Its values at time zero are $\Phi_i(n(w))$, and by part \textup{(i)} of the proposition they differ from $\Phi_i(w)$ by at most $1/i$ in mass. For large $i$, Lemma~\ref{app:return} applies with $\omega=1/i$ and shows that $\mathcal C$ contains a map whose maximal mass is at most $L_i+C_0\mathbf f(\psi_i)$. By the definition of $a$, this gives $a\le L_i+o(1)$. Hence $L_i\to a$.

Let $\mathcal K$ be the set of all varifolds $V=\lim_\ell|\phi_{i_\ell}(v_\ell)|$ with $i_\ell\to\infty$ and $\Mass(\phi_{i_\ell}(v_\ell))\to a$. It is nonempty, because the varifolds of vertices of maximal mass have a convergent subsequence, and it is compact. By Proposition~\ref{app:tightening}\textup{(iii)}, every $V\in\mathcal K$ is stationary with free boundary, and $\|V\|(M)=a$.

\noindent\textbf{Step 2: Local mass decrease.}\enspace
Put $N=N(m)=2\cdot3^m-1$. If some $V\in\mathcal K$ satisfies $\mathcal P_N$, then $V$ is integral by Lemma~\ref{inp:PN-consequences}, and $0<\|V\|(M)=a\le\sup_X\Mass(\Phi_0)<\Lambda$ because $\Phi_0\in\mathcal C$; thus the theorem holds. Suppose, for contradiction, that every $V\in\mathcal K$ fails $\mathcal P_N$. For each $V\in\mathcal K$ there are then a point $p$ and annuli $U_j=\An(p;s_j,r_j)$, $1\le j\le N$, as in \eqref{inp:eq:annular-property}, in none of which $V$ is almost minimizing. Set $\widehat U_j=\An(p;s_j/2,2r_j)$. Since $U_j\Subset\widehat U_j$, Lemma~\ref{inp:localization} implies that $V$ is not almost minimizing with mass steps in any $\widehat U_j$. Thus there is $\zeta_V>0$ such that every relative cycle $T$ with $\dV(|T|,V)<\zeta_V$ admits, for every $\delta_1>0$ and every $j$, a sequence
\begin{equation}\label{inp:eq:local-decrease}
\begin{gathered}
 T=T_0,\ldots,T_q,\quad \spt(T_l-T)\subset\widehat U_j,\\
 \Mass(T_l-T_{l-1})\le\delta_1,\quad
 \Mass(T_l)\le\Mass(T)+\delta_1,\quad
 \Mass(T_q)<\Mass(T)-\zeta_V.
\end{gathered}
\end{equation}
Since $s_j>12r_{j+1}$, the radii $\sigma_j=s_j/2$ and $\rho_j=2r_j$ of the enlarged annuli satisfy \eqref{app:eq:separation}.

The sets $\mathcal O_V=\{W:\dV(W,V)<\zeta_V\}$ are open, and $\mathcal K$ is compact. Hence finitely many of them, say $\mathcal O_1,\ldots,\mathcal O_\mu$, cover $\mathcal K$; let $\zeta>0$ be the smallest of the corresponding numbers $\zeta_V$. We claim that there are $\eta\in(0,a)$ and $i_0$ such that, for $i\ge i_0$, the varifold $|\phi_i(v)|$ of every vertex $v$ with $\Mass(\phi_i(v))\ge a-\eta$ lies in $\mathcal O_1\cup\cdots\cup\mathcal O_\mu$. Otherwise there would be $i_\ell\to\infty$ and vertices $v_\ell$ with $\Mass(\phi_{i_\ell}(v_\ell))\ge a-1/\ell$ whose varifolds lie outside this open set. Since $L_i\to a$, a subsequence of these varifolds would converge to an element of $\mathcal K$ outside the union, which is impossible.

\noindent\textbf{Step 3: Deformation.}\enspace
Fix $i\ge i_0$, and write $\phi=\phi_i$, $j=j_i$ and $\delta_i=\mathbf f(\psi_i)+1/i$, so that $\mathbf f(\phi)\le\delta_i$ and $\delta_i\to0$. Let $W$ be the set of vertices $v\in X(j)_0$ with $\Mass(\phi(v))\ge a-\eta$; for large $i$ we have $a-\eta>m\delta_i$. For each $v\in W$ choose one of the sets $\mathcal O_1,\ldots,\mathcal O_\mu$ containing $|\phi(v)|$. This provides a point $p_v$ and $N$ enlarged annuli about $p_v$, whose radii satisfy \eqref{app:eq:separation}. Consider the graph with vertex set $W$ in which $v\ne v'$ are neighbors if $|v-v'|_\infty\le3^{-j}$. Every vertex has at most $K=3^m-1$ neighbors, and $N=2K+1$. Lemma~\ref{app:selection-lemma} therefore selects for each $v\in W$ one of its annuli, which we call $U_v$, in such a way that $U_v\cap U_{v'}=\varnothing$ for neighbors. By \eqref{inp:eq:local-decrease}, applied with $T=\phi(v)$, the annulus $U_v$ and $\delta_1=\delta_i$, there are sequences $\phi(v)=T^v_0,\ldots,T^v_{q_v}$ satisfying \eqref{app:eq:local-sequence}. Repeating the last term, we may assume that they all have the same length $q$.

Choose $j'>j$ with $3^{j'-j}\ge2q$. Lemma~\ref{app:combination} now gives a discrete homotopy $h_i$ on $X(j')_0$, vanishing on $A$, which starts at $\phi\circ n$ and satisfies
\[
 \mathbf f(h_i)\le(1+2^{m+1})\delta_i,\quad
 \max\Mass(h_i)\le L_i+\kappa_i,\quad
 \kappa_i=2^m(2m+1)\delta_i .
\]
At the final time, a vertex $w$ with $n(w)\in W$ has mass at most $L_i-\zeta+\kappa_i$, and a vertex with $n(w)\notin W$ has mass less than $a-\eta+\kappa_i$. Since $L_i\to a$ and $\kappa_i\to0$, we obtain, for $c=\tfrac12\min\{\zeta,\eta\}$ and all large $i$,
\begin{equation}\label{inp:eq:annular-output}
 \max_w\Mass(h_i(q,w))\le\max\{L_i-\zeta,\ a-\eta\}+\kappa_i\le a-c .
\end{equation}

\noindent\textbf{Step 4: Return to the class $\mathcal C$.}\enspace
Let $\theta_i$ be the discrete homotopy on $X(j')_0$ that first runs through $\psi_i(\cdot,n(w))$ and then through $h_i(\cdot,w)$; the two parts match because $h_i(0,w)=\phi_i(n(w))$. It vanishes on $A$. By Lemma~\ref{app:nearest}\textup{(ii)}, Proposition~\ref{app:tightening}\textup{(ii)} and Step~3, we have $\mathbf f(\theta_i)\to0$ and $\max\Mass(\theta_i)\le a+o(1)$. Its values at time zero are $\Phi_i(n(w))$, and its final values are those of $h_i(q,\cdot)$. As in Step~1, regard $\theta_i$ as a map on the vertex set of $(I\times X)(j'')$ for a large $j''$. Its values at time zero still differ from those of $\Phi_i$ by at most $1/i$ in mass, because the two nearest-vertex maps send a vertex to a vertex of a cell of $X(j_i)$ containing it. For large $i$, Lemma~\ref{app:return} applies with $\omega=1/i$ and produces a map $\widetilde\Phi_i\in\mathcal C$ with
\[
 \sup_X\Mass(\widetilde\Phi_i)\le a-c+C_0\mathbf f(\theta_i)<a .
\]
This contradicts the definition of $a$. Hence some $V\in\mathcal K$ satisfies $\mathcal P_N$, and the theorem is proved.
\end{proof}

\clearpage
\renewcommand{\bibfont}{\small}
\setlength{\bibsep}{3pt}


\begin{thebibliography}{99}
\providecommand{\natexlab}[1]{#1}
\expandafter\ifx\csname urlstyle\endcsname\relax
  \providecommand{\doi}[1]{doi: #1}\else
  \providecommand{\doi}{doi: \begingroup \urlstyle{rm}\Url}\fi

\bibitem[Allard(1972)]{Allard1972}
William~K. Allard.
\newblock On the first variation of a varifold.
\newblock \textit{Annals of Mathematics}, 95\penalty0 (3):\penalty0
  417--491, 1972.

\bibitem[Almgren(1962)]{Almgren1962}
Frederick~J. Almgren, Jr.
\newblock The homotopy groups of the integral cycle groups.
\newblock \textit{Topology}, 1\penalty0 (4):\penalty0 257--299, 1962.

\bibitem[Almgren(1965)]{Almgren1965}
Frederick~J. Almgren, Jr.
\newblock \textit{The Theory of Varifolds}.
\newblock Mimeographed notes, Princeton University, 1965.

\bibitem[Beach(2026)]{Beach2026}
Isabel Beach.
\newblock Short simple orthogonal geodesic chords on a $2$-disk with convex boundary.
\newblock \textit{Journal of Topology and Analysis}, published online
  26 September 2026. \href{https://doi.org/10.1142/S1793525326500500}{DOI 10.1142/S1793525326500500}.

\bibitem[Beach et~al.(2025)Beach, Contreras-Peruyero, Griffin, Rotman, and Searle]{BCGRS2025}
Isabel Beach, Hayde\'e Contreras-Peruyero, Erin Griffin, Regina Rotman, and Catherine Searle.
\newblock Lengths of orthogonal geodesic chords on Riemannian manifolds.
\newblock arXiv:2509.19620v1, 2025.

\bibitem[Bellettini et~al.(2026)Bellettini, Li, Parise, and Sarnataro]{BLPS2026}
Costante Bellettini, Martin Man-chun Li, Davide Parise, and Lorenzo Sarnataro.
\newblock Regularity and compactness theory of stable free-boundary minimal
  hypersurfaces.
\newblock arXiv:2610.09596v1, 2026.

\bibitem[Cheeger and Colding(1996)]{CC1996}
Jeff Cheeger and Tobias~H. Colding.
\newblock Lower bounds on Ricci curvature and the almost rigidity of
  warped products.
\newblock \textit{Annals of Mathematics},
  144\penalty0 (1):\penalty0 189--237, 1996.

\bibitem[Cheeger and Naber(2015)]{CN2015}
Jeff Cheeger and Aaron Naber.
\newblock Regularity of Einstein manifolds and the codimension 4 conjecture.
\newblock \textit{Annals of Mathematics}, 182\penalty0 (3):\penalty0
  1093--1165, 2015.
\newblock Numbering follows the published version.

\bibitem[Chu and Li(2025)]{CL2025}
Adrian Chun-Pong Chu and Yangyang Li.
\newblock Strong multiplicity one theorems and homological min-max theory.
\newblock \textit{Advances in Mathematics}, 482\penalty0 (Part B):\penalty0
  110623, 2025.

\bibitem[Colding et~al.(2018)Colding, Gabai, and Ketover]{CGK2018}
Tobias~H. Colding, David Gabai, and Daniel Ketover.
\newblock On the classification of {Heegaard} splittings.
\newblock \textit{Duke Mathematical Journal}, 167\penalty0 (15):\penalty0
  2833--2856, 2018.

\bibitem[De~Lellis and Ramic(2018)]{DLR2018}
Camillo De~Lellis and Jusuf Ramic.
\newblock Min-max theory for minimal hypersurfaces with boundary.
\newblock \textit{Annales de l'Institut Fourier}, 68\penalty0 (5):\penalty0
  1909--1986, 2018.

\bibitem[Evans and Gariepy(1992)]{EG1992}
Lawrence~C. Evans and Ronald~F. Gariepy.
\newblock \textit{Measure Theory and Fine Properties of Functions}.
\newblock Studies in Advanced Mathematics. CRC Press, Boca Raton, FL, 1992.

\bibitem[Federer(1969)]{Federer1969}
Herbert Federer.
\newblock \textit{Geometric Measure Theory}, volume 153 of \textit{Grundlehren der
  mathematischen Wissenschaften}.
\newblock Springer, 1969.

\bibitem[Fu and Zhu(2026)]{FZ2026}
Wenjie Fu and Zhifei Zhu.
\newblock Homological filling and minimal varifolds in four-dimensional
  Einstein manifolds.
\newblock arXiv:2512.14016v2, revised 6 March 2026.
\newblock Numbering follows this revised version.

\bibitem[Gallot(1988)]{Gallot1988}
Sylvestre Gallot.
\newblock In\'egalit\'es isop\'erim\'etriques et analytiques sur les
  vari\'et\'es riemanniennes.
\newblock In \textit{On the Geometry of Differentiable Manifolds},
  \textit{Ast\'erisque}, 163--164:\penalty0 31--91, 1988.

\bibitem[Glynn-Adey and Liokumovich(2017)]{GAL2017}
Parker Glynn-Adey and Yevgeny Liokumovich.
\newblock Width, Ricci curvature, and minimal hypersurfaces.
\newblock \textit{Journal of Differential Geometry},
  105\penalty0 (1):\penalty0 33--54, 2017.
\newblock arXiv:1408.3656v4.

\bibitem[Gromov(1983)]{Gromov1983}
Mikhael Gromov.
\newblock Filling {Riemannian} manifolds.
\newblock \textit{Journal of Differential Geometry}, 18\penalty0 (1):\penalty0
  1--147, 1983.

\bibitem[Gr{\"u}ter and Jost(1986)]{GJDisks1986}
Michael Gr\"uter and J\"urgen Jost.
\newblock On embedded minimal disks in convex bodies.
\newblock \textit{Annales de l'Institut Henri Poincar\'e C, Analyse non
  lin\'eaire}, 3\penalty0 (5):\penalty0 345--390, 1986.

\bibitem[Guang et~al.(2021)Guang, Li, Wang, and Zhou]{GLWZ2021}
Qiang Guang, Martin Man-chun Li, Zhichao Wang, and Xin Zhou.
\newblock Min-max theory for free boundary minimal hypersurfaces {II}: General
  {Morse} index bounds and applications.
\newblock \textit{Mathematische Annalen}, 379:\penalty0 1395--1424, 2021.

\bibitem[Guth(2007)]{Guth2007}
Larry Guth.
\newblock The width-volume inequality.
\newblock \textit{Geometric and Functional Analysis}, 17\penalty0 (4):\penalty0
  1139--1179, 2007.

\bibitem[Li(2015)]{Li2015}
Martin Man-chun Li.
\newblock A general existence theorem for embedded minimal surfaces with free
  boundary.
\newblock \textit{Communications on Pure and Applied Mathematics},
  68\penalty0 (2):\penalty0 286--331, 2015.

\bibitem[Li and Zhou(2021)]{LZ2021}
Martin Man-chun Li and Xin Zhou.
\newblock Min-max theory for free boundary minimal hypersurfaces {I}:
  Regularity theory.
\newblock \textit{Journal of Differential Geometry}, 118\penalty0 (3):\penalty0
  487--553, 2021.

\bibitem[Liokumovich et~al.(2018)Liokumovich, Marques, and Neves]{LMN2018}
Yevgeny Liokumovich, Fernando~C. Marques, and Andr{\'e} Neves.
\newblock Weyl law for the volume spectrum.
\newblock \textit{Annals of Mathematics}, 187\penalty0 (3):\penalty0
  933--961, 2018.

\bibitem[Liokumovich et~al.(2015)Liokumovich, Nabutovsky, and Rotman]{LNR2015}
Yevgeny Liokumovich, Alexander Nabutovsky, and Regina Rotman.
\newblock Contracting the boundary of a {Riemannian} $2$-disc.
\newblock \textit{Geometric and Functional Analysis}, 25\penalty0 (5):\penalty0
  1543--1574, 2015.

\bibitem[Marchese and Stuvard(2018)]{MS2018}
Andrea Marchese and Salvatore Stuvard.
\newblock On the structure of flat chains modulo $p$.
\newblock \textit{Advances in Calculus of Variations},
  11\penalty0 (3):\penalty0 309--323, 2018.

\bibitem[Marques and Neves(2014)]{MN2014}
Fernando~C. Marques and Andr{\'e} Neves.
\newblock Min-max theory and the Willmore conjecture.
\newblock \textit{Annals of Mathematics},
  179\penalty0 (2):\penalty0 683--782, 2014.

\bibitem[Nabutovsky and Rotman(2003)]{NR2003Sphere}
Alexander Nabutovsky and Regina Rotman.
\newblock The area of a minimal embedded $2$-sphere in a manifold
  diffeomorphic to $S^3$.
\newblock \textit{International Mathematics Research Notices},
  2003\penalty0 (39):\penalty0 2121--2129, 2003.

\bibitem[Nabutovsky and Rotman(2004)]{NR2004}
Alexander Nabutovsky and Regina Rotman.
\newblock Volume, diameter and the minimal mass of a stationary $1$-cycle.
\newblock \textit{Geometric and Functional Analysis}, 14:\penalty0 748--790,
  2004.

\bibitem[Nabutovsky and Rotman(2006)]{NR2006}
Alexander Nabutovsky and Regina Rotman.
\newblock Curvature-free upper bounds for the smallest area of a minimal
  surface.
\newblock \textit{Geometric and Functional Analysis}, 16\penalty0 (2):\penalty0
  453--475, 2006.
\newblock Numbering follows the manuscript of 20 June 2005.

\bibitem[Nabutovsky et~al.(2021)Nabutovsky, Rotman, and Sabourau]{NRS2020}
Alexander Nabutovsky, Regina Rotman, and St{\'e}phane Sabourau.
\newblock Sweepouts of closed {Riemannian} manifolds.
\newblock \textit{Geometric and Functional Analysis}, 31\penalty0 (3):\penalty0
  721--766, 2021.

\bibitem[Pitts(1981)]{Pitts1981}
Jon~T. Pitts.
\newblock \textit{Existence and Regularity of Minimal Surfaces on {Riemannian}
  Manifolds}, volume~27 of \textit{Mathematical Notes}.
\newblock Princeton University Press, Princeton, NJ, 1981.

\bibitem[Rotman(2011)]{Rotman2011}
Regina Rotman.
\newblock Flowers on {Riemannian} manifolds.
\newblock \textit{Mathematische Zeitschrift}, 269:\penalty0 543--554, 2011.

\bibitem[Sarnataro et~al.(2026)Sarnataro, Stryker, Wang, and Zhou]{SSWZ2026}
Lorenzo Sarnataro, Douglas Stryker, Zhichao Wang, and Xin Zhou.
\newblock Existence of free boundary minimal disks in convex regions.
\newblock arXiv:2606.02396, 2026.

\bibitem[Simon(2014)]{Simon2014}
Leon Simon.
\newblock \textit{Introduction to Geometric Measure Theory}.
\newblock Lecture notes, Stanford University, 2014.

\bibitem[Struwe(1984)]{Struwe1984}
Michael Struwe.
\newblock On a free boundary problem for minimal surfaces.
\newblock \textit{Inventiones Mathematicae}, 75\penalty0 (3):\penalty0
  547--560, 1984.

\bibitem[Sun et~al.(2024)Sun, Wang, and Zhou]{SWZ2024}
Ao Sun, Zhichao Wang, and Xin Zhou.
\newblock Multiplicity one for min--max theory in compact manifolds with
  boundary and its applications.
\newblock \textit{Calculus of Variations and Partial Differential Equations},
  63, 2024.

\bibitem[Whitney(1957)]{Whitney1957}
Hassler Whitney.
\newblock \textit{Geometric Integration Theory}.
\newblock Princeton University Press, 1957.

\bibitem[Wu and Zhu(2020)]{WZ2020}
Nan Wu and Zhifei Zhu.
\newblock An upper bound for the smallest area of a minimal surface
  in manifolds of dimension four.
\newblock \textit{Journal of Geometric Analysis},
  30:\penalty0 573--600, 2020.
\newblock Numbering follows arXiv:1702.07034v2.

\bibitem[Zhang and Zhu(2026)]{PartI}
Tengyao Zhang and Zhifei Zhu.
\newblock Free-boundary minimal varifolds I: zero boundary mass.
\newblock \href{https://arxiv.org/abs/2611.30322}{arXiv:2611.30322}, 2026.

\bibitem[Zhou(2016)]{Zhou2016}
Xin Zhou.
\newblock On the free boundary min-max geodesics.
\newblock \textit{International Mathematics Research Notices},
  2016\penalty0 (5):\penalty0 1447--1466, 2016.

\bibitem[Zhu(2023)]{JZ2023}
Jonathan~J. Zhu.
\newblock Widths of balls and free boundary minimal submanifolds.
\newblock \textit{Advanced Nonlinear Studies}, 23\penalty0 (1):\penalty0
  20220044, 2023.

\end{thebibliography}
\end{document}